\documentclass[12pt]{article}
\usepackage{amsthm,amsfonts,amssymb,amscd}

\begin{document}


\begin{center} \large \bf Factoriality and birational rigidity \\
of two families of singular quartic three-folds
\end{center}\vspace{0.5cm}

\centerline{A.V.Pukhlikov}\vspace{0.5cm}

\parshape=1
3cm 11cm \noindent {\small \quad\quad\quad \quad\quad\quad\quad
\quad\quad\quad {\bf }\newline In this paper we study two families
of three-dimensional quartics in the complex projective space
${\mathbb P}^4$: hypersurfaces with a unique quadratic singularity
of rank 3, which is resolved by two blowups, and hypersurfaces
with two quadratic singularities of rank 3 and 4. Both families
have codimension 3 in the natural parameter space. For a Zariski
general quartic in each of these families we prove factoriality
and birational rigidity and describe its group of birational
self-maps.

Bibliography: 15 items.} \vspace{1cm}


\vspace{1cm}

\section*{Introduction}

{\bf 0.1. The main result.} The present paper is a follow up of
\cite{Pukh2025a}. We prove factoriality and birational rigidity
and describe the group of birational automorphisms of two classes
of three-dimensional quartics in ${\mathbb P}^4$ (the ground field
is the field of complex numbers ${\mathbb C}$):

(1) the quartics $V\subset {\mathbb P}^4$ with a unique singular
point $o_1\in V$, which is a quadratic singularity of rank 3, the
blowup  $\varphi_1\colon V_1\to V$ of which yields a variety $V_1$
with a unique singularity, a non-degenerate quadratic singularity
$o_2\in V_1$ (so that the blowup $\varphi_2\colon V_2\to V_1$ of
the point $o_2$ gives a non-singular variety $V_2$), --- in this
paper we will call such quartics the quartics {\it of the first
type},

(2) the quartics $V\subset {\mathbb P}^4$ with two quadratic
singular points $o_1$, $o_2$, where the rank of the point $o_1$ is
3, and the rank of the point $o_2$ is 4, and the blowup of the two
points $o_1$, $o_2$ gives a non-singular variety, we will call
such hypersurfaces the quartics {\it of the second type}.

In both cases we will assume that the quartic $V$ satisfies
certain conditions of general position, see Subsection 0.2. Both
families have codimension 3 in the (projectively) 69-dimensional
space of three-dimensional quartics. Thus the present paper can be
considered as the next step (after \cite{Pukh2025a}) in the proof
of the old conjecture \cite{CoMe,AbbanKaloghiros2016} that
factorial three-dimensional quartics with isolated singularities
not worse than $cA_1$ are birationally rigid; or as the next step
in the description of the {\it birational geography} of
three-dimensional quartics, see Subsection 0.5
below.\vspace{0.3cm}


{\bf 0.2. Quadratic singularities.} Our notations are as close as
possible to the notations of \cite{Pukh2025a}. Let ${\mathbb
A}^4\subset {\mathbb P}^4$ be an affine chart with the coordinates
$(z_1,z_2,z_3,z_4)$ and
$$
f(z_*)=q_2(z_*)+q_3(z_*)+q_4(z_*)
$$
a polynomial of degree 4, where $q_i(z_*)$ are homogeneous of
degree $i\in\{2,3,4\}$, where $q_2=z_1^2+z_2^2+z_3^2$ and the set
of common zeros of the polynomials $q_2$, $q_3$, $q_4$ on
${\mathbb P}^3_{(z_1:z_2:z_3:z_4)}$ consists of 24 distinct
points. Let $\pi_1\colon X_1\to {\mathbb P}^4$ be the blow up of
the point $o_1=(0,0,0,0)$, denote by the symbol ${\mathbb E}_1$
its exceptional divisor, ${\mathbb E}_1\cong {\mathbb
P}^3_{(z_1:z_2:z_3:z_4)}$. The symbol $V$ stands for the
hypersurface $\{f=0\}\subset {\mathbb P}^4$, and the symbol $V_1$
--- for its strict transform on $X_1$, so that $\varphi_1\colon
V_1\to V$ is the blow-up of the point $o_1$ (now on $V$). Let
$Q_1=V_1\cap {\mathbb E}_1$ be the exceptional quadric. The only
candidate for a singular point of the variety $V_1$ on $Q_1$ is
the vertex of the quadratic cone $Q_1$, i.e. the point
$$
o_2=(0:0:0:1).
$$
We define the quartic of the first type, assuming that $o_2\in
V_1$ is a non-degenerate (nodal) quadratic point: in terms of the
polynomial $f$ it means that
$$
q_3(0,0,0,1)=0,
$$
however,
$$
4q_4(0,0,0,1)-\sum_{i=1}^3\left(\frac{\partial^3 q_3}{\partial
z_i\partial z_4^2}(0,0,0,1)\right)^2\neq 0
$$
(the equality of the left hand side to zero means precisely that
the point $o_2$ is a quadratic singularity of rank 3). Therefore,
the singular point $o_1\in V$ is resolved by two blow-ups: let
$\pi_2\colon X_2\to X_1$ be the blow-up of the point $o_2$ with
the exceptional divisor ${\mathbb E}_2\cong {\mathbb P}^3$ and
$V_2\subset X_2$ the strict transform of the hypersurface $V_1$,
that is, $\varphi_2\colon V_2\to V_1$ is the blow-up of the point
$o_2$ on $V_1$, then $V_2\cap {\mathbb E}_2=Q_2$ is a
non-degenerate quadric. We assume that
$$
q_4(0,0,0,1)\neq 0
$$
(geometrically it means that the line $[o_1,o_2]\subset {\mathbb
P}^4$, passing through the point $o_1$ in the direction $o_2$, is
not contained in $V$). We assume that $\{o_1\}=\mathop{\rm Sing}
V$. In order to reduce the number of cases to be considered, we
will assume that the quartic $V$ satisfies some additional
conditions of general position, which will be stated in \S\S 3,5.

The quartics of the second type are defined by the conditions
$q_3(0,0,0,1)\neq 0$ (so that $V_1$ is non-singular at the points
of the quadric $Q_1$), however, $V$ has a non-degenerate quadratic
singularity $o_2\neq o_1$, and moreover,
$$
\mathop{\rm Sing} V=\{o_1,o_2\}.
$$
In addition we assume that the line $[o_1,o_2]\subset {\mathbb
P}^4$, containing these points, is not contained in $V$. Apart
from that, certain conditions of general position, simplifying our
arguments, are stated in \S\S 2,4.

In order to make our notations uniform, in the case of a quartic
of the second type let $\pi_2\colon X_2\to X_1$ be the blow-up of
the point $o_2$ (which now is not infinitely near to $o_1$),
${\mathbb E}_2\cong {\mathbb P}^3$ the exceptional divisor,
$V_2\subset X_2$ the strict transform and $Q_2=V_2\cap {\mathbb
E}_2$ a non-degenerate quadric. The second blow-up $\pi_2$ does
not affect the cone $Q_1$ (as $o_2\not\in Q_1$) and one may write
$Q_1\subset X_2$; obviously, $Q_1\cap Q_2=\emptyset$.

{\bf Remark 0.1.} For quartics of the first type the condition
$[o_1,o_2]\not\subset V$ follows from the assumption that on $V$
there are precisely 24 distinct lines, passing through the point
$o_1$: if $q_4(0,0,0,1)=0$, then the vertex of the cone $Q_1$ (the
point $o_2$) comes into the cycle of the scheme-theoretic
intersection of the surfaces
$$
\{q_2=0\},\quad \{q_3=0\},\quad \{q_4=0\}
$$
in ${\mathbb E}_1\cong {\mathbb P}_{(z_1:z_2:z_3:z_4)}$ with
multiplicity $\geqslant 2$, so that on $V$ there are at most 23
distinct lines, passing through the point $o_1$. For quartics of
the second type the condition $[o_1,o_2]\not\subset V$ is assumed
to hold in order to reduce the number of cases to inspect in the
proof of the plane section lemma, see \S 4.\vspace{0.3cm}


{\bf 0.3. Factoriality and birational rigidity.} We have the
equality
$$
\mathop{\rm Pic} X_2={\mathbb Z} H_{{\mathbb P}}\oplus {\mathbb Z}
{\mathbb E}_1 \oplus {\mathbb Z} {\mathbb E}_2
$$
(which holds for quartics of both types), where we omit the
symbols of pull back $\pi_1^*$, $\pi_2^*$, and $H_{{\mathbb P}}$
is the class of a hyperplane in ${\mathbb P}^4$. Obviously,
$$
{\mathbb Z} H\oplus {\mathbb Z} Q_1\oplus {\mathbb Z} Q_2\subset
\mathop{\rm Pic} V_2
$$
(here $H$ is the class of a hyperplane section of $V$).

{\bf Theorem 0.1.} {\it The following equality holds:
$$
\mathop{\rm Pic} V_2={\mathbb Z} H\oplus {\mathbb Z} Q_1\oplus
{\mathbb Z} Q_2,
$$
so that $V\subset {\mathbb P}^4$ is a factorial hypersurface with}
$\mathop{\rm Cl} V=\mathop{\rm Pic} V={\mathbb Z} H$.

Let $\Sigma$ be a mobile linear system on $V$ (that is,
$\mathop{\rm codim} \mathop{\rm Bs} \Sigma\geqslant 2$, the system
$\Sigma$ has no fixed components). By Theorem 0.1, there is
$n[\Sigma]\geqslant 1$, such that $\Sigma\subset |n[\Sigma]\cdot
H|$. We say that $\Sigma$ {\it has no maximal singularities}, if
the pair $\left(V, \frac{1}{n[\Sigma]} D\right)$ is canonical,
that is, for every exceptional divisor $E$ over $V$ the inequality
$$
\mathop{\rm ord}\nolimits_E D\leqslant n[\Sigma]\cdot a(E)
$$
holds, where $D\in\Sigma$ is a general divisor and $a(E)$ is the
discrepancy.

{\bf Theorem 0.2.} {\it The quartic $V$ is a birationally rigid
variety: for every mobile linear system $\Sigma$ there is a
birational automorphism $\chi\in\mathop{\rm Bir} V$, such that the
linear system $(\chi^{-1})_*\Sigma$ has no maximal singularities.
On $V$ there are no structures of a conic bundle and fibrations
into del Pezzo surfaces, and every birational map
$V\dashrightarrow V'$ onto a three-dimensional Fano variety with
${\mathbb Q}$-factorial terminal singularities is an isomorphism.}

(It is well known that the second claim of Theorem 0.2 follows
from the first one.)

The third main result of the present paper is an explicit
description of the group of birational automorphisms $\mathop{\rm
Bir} V$. In contrast to the first two main results, its statement
is not the same for the quartics of the first and second type.
This is because the second singularity of a quartic of the first
type is an infinitely near point and does not generate any
additional birational automorphisms.

For a quartic $V$ of the first type let $\tau_0$ be the Galois
involution, generated by the projection $V\dashrightarrow {\mathbb
P}^3$ from the point $o_1$ and $\tau_1$, \dots, $\tau_{24}$ are
the birational involutions, linked to the 24 lines on $V$, passing
through the point $o_1$, see \cite{Pukh89c,Pukh2025a}. The
subgroup in $\mathop{\rm Bir} V$, generated by these 25
involutions, will be denoted by the symbol $B(V)$. For a quartic
$V$ of the second type let $\tau_{1,0},\dots, \tau_{1,24}$ have
the same meaning as the 25 involutions, listed above, $\tau_{2,0}$
the Galois involution, generated by the projection from the point
$o_2$, and $\tau_{2,i}$, where $i\in\{1,\dots, L(o_2)\}$, are the
birational involutions, linked to the lines on $V$, passing
through the point $o_2$. We know that $L(o_2)\leqslant 24$.

{\bf Proposition 0.1.} {\it For a general quartic of the second
type the equality} $L(o_2)= 24$ {\it holds.}

We denote by the symbol $B(V)$ the subgroup in $\mathop{\rm Bir}
V$, generated by these $26+L(o_2)$ involutions.

{\bf Theorem 0.3.} {\it The subgroup $B(V)$ is normal in
$\mathop{\rm Bir} V$, and there is the exact sequence
$$
1\to B(V)\to \mathop{\rm Bir} V \to \mathop{\rm Aut} V\to 1,
$$
where the subgroup $\mathop{\rm Aut} V\subset \mathop{\rm Bir} V$
is finite, for a general (in its family) quartic it is trivial.
The involutions, generating the subgroup $B(V)$, form a system of
free generators:
$$
B(V)=\mathop{*}\limits_{i=0}^{24}\langle \tau_i\rangle
$$
for quartics of the first type and
$$
B(V)=\mathop{*}\limits_{i=0}^{24}\langle \tau_{1,i}\rangle *\,\,
\mathop{*}\limits_{i=0}^{L(o_2)}\langle \tau_{2,i}\rangle
$$
for quartics of the second type.}

{\bf Remark 0.2.} In the semi-direct product for the group
$\mathop{\rm Bir} V$ the finite subgroup $\mathop{\rm Aut} V$ acts
on $B(V)$ by conjugation, preserving the point $o_1$ and permuting
the lines, passing through the point $o_1$, for quartics of the
first type, and preserving both points $o_1$, $o_2$ and permuting
the lines passing through these points (inside each of the two
groups of lines). Thus in the first family of quartics the
subgroup $B(V)$ is freely generated by 25 involutions and in the
second family for a general quartic it is freely generated by 50
involutions.\vspace{0.3cm}


{\bf 0.4. The structure of the paper.} In \S 1 we prove that
quartics of the first and second type are factorial (Theorem 0.1).
In \S 2 we start the proof of Theorems 0.2 and 0.3: we recall the
concept of a maximal singularity of a mobile linear system of
divisors on the quartic $V$ and carry out the standard preparatory
work (list the possible types of maximal singularities and exclude
some of them). The claims of Theorems 0.2 and 0.3 for quartics of
the first type are reduced to two facts: the plane section lemma
(Proposition 2.2) and the claim that the centre of a maximal
singularity on $V_1$ can not be a generator of the cone $Q_1$
(Proposition 2.4). For quartics of the second type it is
sufficient to show the plane section lemma and a few close claims,
some of which are shown in \S 2, and the others are considered in
\S 4.

In \S 3 we show the plane section lemma for quartics of the first
type. We go through all options for quartics of general position
(in that family) and a few additional examples. Apart from that,
we show some auxiliary claims that are used later when we exclude
generators of the cone $Q_1$ in \S 5 (this is the hardest part of
the proof of Theorems 0.2 and 0.3).

In \S 4 we show the plane section lemma for quartics of the second
type. Here we also prove Proposition 2.5. Its proof requires the
same technique as that of the plane section lemma.

Finally, in \S 5 we show the key intermediate claim in the proof
of the birational rigidity of quartics of the first type: a
generator $R$ of the cone $Q_1$ can not be the centre of a maximal
singularity. We use both the quadratic method (the technique of
counting multiplicities for the self-intersection of the system
$\Sigma$), and the linear one (similar to the proof of the plane
section lemma).\vspace{0.3cm}


{\bf 0.5. Birational geography of three-dimensional quartics.}
Irreducible quartic hypersurfaces in ${\mathbb P}^4$ are
parameterized by points of a Zariski open subset in the projective
space ${\mathbb P}^{69}$, that is, by points of a quasi-projective
69-dimensional variety. In \cite[Subsection 0.4]{Pukh2025a} the
problem of studying the ``birational geography'' of
three-dimensional quartics was formulated: to break the natural
parameter space into disjoint quasi-projective subvarieties,
consisting of quartics with a similar birational type. The natural
``measure of size'' of these subvarieties is their dimension (or
the codimension in ${\mathbb P}^{69}$). The seminal paper
\cite{IM} that started the theory of birational rigidity described
the birational type of non-singular quartics, reducing the
subsequent work to investigating the 68-dimensional subset of
singular quartics, that is, to the subset of codimension 1. It
follows from \cite{Pukh89c} that the birational rigidity still
holds in codimension 1, too (although the group of birational
automorphisms becomes very large). It follows from \cite{Me04}
(taking into account \cite{Shr08b,Ch2010,Shramov2007}) and
\cite{Pukh2025a} that the codimension of the complement to the set
of factorial birationally rigid three-dimensional quartics is at
least 3. Non-rigid but ``almost rigid'' quartics were studies in
\cite{CoMe}, see also \cite{AbbanKaloghiros2016}. For the families
and types of singularities that are already accessible for the
modern techniques (however, often very hard to study), see.
\cite[Subsection 0.4]{Pukh2025a}. Note also the paper
\cite{Krylovetal2024}, which it is natural to consider as a part
of the birational geography of Fano double spaces of dimension 3,
branched over a sextic surface in ${\mathbb P}^3$.




\section{Factoriality}

In this section we prove Theorem 0.1.\vspace{0.3cm}

{\bf 1.1. Preliminary remarks.} We use the notations of
Subsections 0.2, 0.3, considering in detail quartics of the first
type. Note, in the first place, that the proof of factoriality,
presented in \cite[Remark 1.2]{CoMe}, works for a quartics of the
first type $V$, too, as by construction the homogeneous polynomial
$q_4(z_*)$ is arbitrary: the section of $V$ by the hypersurface
``at infinity'' is precisely $\{q_4=0\}$. Therefore, for a general
quartic of the first type its section ``at infinity'' has Picard
number 1. However, we give a different proof of factoriality that
does not use this fact.

Note, first of all, that on $X_2$ the class of the hypersurface
$V_2$ is $4H_{\mathbb P}-2{\mathbb E}_1-2{\mathbb E}_2$. This
class is nef and big, but not ample, so that the Lefschetz theorem
does not apply (in contrast to the quartics considered in
\cite{Pukh2025a}). Furthermore, it is easy to see that all
cohomology spaces
$$
H^1({\mathbb P}^4,{\cal O}_{{\mathbb P}^4}),\quad H^1(X_1,{\cal O}_{X_1}),
\quad H^1(X_2,{\cal O}_{X_2}),
$$
$$
H^2({\mathbb P}^4,{\cal O}_{{\mathbb P}^4}),\quad H^2(X_1,{\cal O}_{X_1}),
\quad H^2(X_2,{\cal O}_{X_2}),
$$
and the six similar spaces for $V$, $V_1$, $V_2$ are trivial, so
that (by the exponential sequences for $X_2$ and $V_2$) there are
the natural isomorphisms
$$
\mathop{\rm Pic} X_2\otimes {\mathbb C}\cong H^1(X_2,\Omega_{X_2})
$$
and
$$
\mathop{\rm Pic} V_2\otimes {\mathbb C}\cong H^1(V_2,\Omega_{V_2}),
$$
and moreover, the restriction of 1-forms onto the non-singular
hypersurface $V_2$ determines the inclusion
$$
\rho\colon H^1(X_2,\Omega_{X_2})\to H^1(V_2,\Omega_{V_2})
$$
and we have to show that $\rho$ is an isomorphism.\vspace{0.3cm}


{\bf 1.2. Exact sequences.} Consider the exact sequence
$$
0\to {\cal N}^*_{V_2/X_2}={\cal O}_{X_2}(-V_2)|_{V_2}\to
\Omega_{X_2}|_{V_2}\to \Omega_{V_2}\to 0.
$$
From the obvious equalities $H^0(V_2,\Omega_{V_2})=0$ (so that
$H^0(V_2,\Omega_{X_2}|_{V_2})=0$, too) and
$$
H^1(V_2,{\cal O}_{V_2}(-V_2))=H^2(V_2,{\cal O}_{V_2}(-V_2))=0
$$
(since $4H-2Q_1-2Q_2$  is a nef and big divisor on $V_2$) we get
the natural isomorphism
$$
\gamma\colon H^1(V_2,\Omega_{X_2}|_{V_2})\to H^1(V_2,\Omega_{V_2}).
$$
Now consider the exact sequence
$$
0\to{\cal O}_{X_2}(-V_2)\to{\cal O}_{{X_2}}\to{\cal O}_{{V_2}}\to 0,
$$
tensored by $\Omega_{X_2}$. We get the exact sequence
$$
0\to H^1(X_2,\Omega_{X_2}(-V_2))\to H^1(X_2,\Omega_{X_2})\to
H^1(V_2,\Omega_{X_2}|_{V_2})\to
$$
$$
\to H^2(X_2,\Omega_{X_2}(-V_2)) \to H^2(X_2,\Omega_{X_2})=0.
$$
Since the isomorphism $\gamma$ is a natural one, we conclude that
in that exact sequence the map
$$
H^1(X_2,\Omega_{X_2})\to H^1(V_2,\Omega_{X_2}|_{V_2})
$$
is injective (it is the inclusion $\rho$, twisted by the
isomorphism $\gamma$), so that
$$
H^1(X_2,\Omega_{X_2}(-V_2))=0.
$$
If we also have $H^2(X_2,\Omega_{X_2}(-V_2))=0$, then we obtain
the required claim. Note that since $4H_{\mathbb P}-2{\mathbb
E}_1-2{\mathbb E}_2$ is not ample, although it is nef and big, the
Akizuki-Nakano theorem \cite{AkizukiNakano1954} does not apply,
see also \cite{EsnaultViehweg1992}.

Now consider a general quadric $S\subset{\mathbb P}^4$, such that
$S\ni o_1$ and the strict transform $S_1\subset X_1$ contains the
point $o_2$. Denote the strict transform of $S_1$ on $X_2$ by the
symbol $S_2$. Obviously, the hypersurfaces $S,S_1,S_2$ are
non-singular and the class of divisor $S_2$ on $X_2$ is
$$
2H_{\mathbb P}-{\mathbb E}_1-{\mathbb E}_2,
$$
so that $V_2\sim 2S_2$. Since the original quadric $S$ is
non-singular, we have
$$
\mathop{\rm Pic} S_2=
{\mathbb Z}H_S\oplus{\mathbb Z}P_1\oplus{\mathbb Z}P_2,
$$
where $H_S$ is the class of a hyperplane section and $P_1={\mathbb
E}_1|_{S_1}$, $P_2={\mathbb E}_2|_{S_2}$ (recall that the pull
back of a divisor is denoted by the same symbol). Repeating the
arguments, performed above for the pair $X_2$, $V_2$, for the pair
$X_2$, $S_2$, we obtain the exact sequence
$$
0\to H^1(X_2,\Omega_{X_2}(-S_2))\to H^1(X_2,\Omega_{X_2})\to
H^1(S_2,\Omega_{X_2}|_{S_2})\to H^2(X_2,\Omega_{X_2}(-S_2)) \to 0.
$$
However, the central map in that sequence is induced by the
restriction map $\mathop{\rm Pic} X_2\to\mathop{\rm Pic} S_2$,
which is an isomorphism since the quadric $S$ is smooth, so we
conclude that
$$
H^i(X_2,\Omega_{X_2}(-S_2))=0
$$
for $i=1,2$. Taking these equalities into account, from the short
exact sequence
$$
0\to{\cal O}_{X_2}(-S_2)\to{\cal O}_{X_2}\to{\cal O}_{S_2}\to 0,
$$
tensored by $\Omega_{X_2}(-S_2)$, we get the isomorphism
$$
H^1(S_2,\Omega_{X_2}(-S_2)|_{S_2})\cong H^2(X_2,\Omega_{X_2}(-V_2)).
$$

{\bf Lemma 1.1.} {\it The space $H^1(S_2,\Omega_{X_2}(-S_2)|_{S_2})$
is trivial.}

{\bf Proof} is given in Subsection 1.3.

Lemma 1.1 obviously implies the claim of Theorem 0.1 for quartics
of the first type.\vspace{0.3cm}


{\bf 1.3. Proof of Lemma 1.1.} Our arguments are similar to the
arguments of the previous subsection. Tensoring the exact sequence
$$
0\to {\cal O}_{S_2}(-S_2)\to \Omega_{X_2}|_{S_2}\to \Omega_{S_2}\to 0
$$
by ${\cal O}_{S_2}(-S_2)$ and taking into account that the class
$$
S_2|_{S_2}=2H_S-P_1-P_2
$$
is nef and big, we obtain the equalities
$$
H^1(S_2,{\cal O}_{S_2}(-2S_2))=H^2(S_2,{\cal O}_{S_2}(-2S_2))=0,
$$
so that
$$
H^1(S_2,\Omega_{X_2}(-S_2)|_{S_2})=H^1(S_2,\Omega_{S_2}(-S_2))
$$
and we have to show that the space $H^1(S_2,\Omega_{S_2}(-S_2))$
is trivial. Consider the free linear system $|2H_S-P_1-P_2|$, let
$W$ be a general surface in that system. Obviously, its image in
${\mathbb P}^4$ is a non-singular del Pezzo surface of degree 4,
and $W$ is its blow-up first at the point $o_1$, and then (of its
strict transform) at the point $o_2$. The rank of the group
$\mathop{\rm Pic} W$ is 8, and the restriction map
$$
\rho_W\colon \mathop{\rm Pic}S_2\to \mathop{\rm Pic} W
$$
is injective. For the surface $W$ the equalities $H^i(W,{\cal
O}_W)=0$ hold for $i=1,2$, so that $\mathop{\rm Pic}
W\otimes{\mathbb C}=H^1(W,\Omega_W)$. The injective map
$$
H^1(S_2,\Omega_{S_2})\to H^1(W,\Omega_W),
$$
induced by the restriction of the forms, we denote by the same
symbol $\rho_W$. Since the divisor $(2H_S-P_1-P_2)|_W$ is nef and
big, we get the equality $H^1(W,{\cal O}_W(-S_2))=0$, so that the
exact sequence
$$
0\to {\cal O}_{W}(-S_2)\to \Omega_{S_2}|_{W}\to \Omega_{W}\to 0
$$
implies injectivity of the map
$$
\beta\colon H^1(W,\Omega_{S_2}|_{W})\to H^1(W,\Omega_{W})
$$
(and also the equality $H^0(W,\Omega_{S_2}|_{W})=0$). Now let us
tensor the exact sequence
$$
0\to {\cal O}_{S_2}(-S_2)\to {\cal O}_{S_2}\to {\cal O}_{W}\to 0
$$
by $\Omega_{S_2}$ and obtain an exact sequence of linear maps
$$
0\to H^1(S_2,\Omega_{S_2}(-S_2))\to H^1(S_2,\Omega_{S_2})
\stackrel{\alpha}{\to} H^1(W,\Omega_{S_2}|_{W}).
$$
Since the composition $\beta\circ\alpha=\rho_W$ is injective, we
conclude from here that the space $H^1(S_2,\Omega_{S_2}(-S_2))=0$,
as we claimed. Q.E.D. for Lemma 1.1. \vspace{0.3cm}


{\bf 1.4. Quartics of the second type.} Now let us make sure that
the arguments above work for quartics of the second type. Since
the notations for both cases are almost identical, no changes in
the formulas are needed. As we will show below (see Subsection
4.3), in the case of a second quadratic singularity $o_2\neq o_1$
the form $q_4$ can be arbitrary (there are some relations between
$q_3$ and $q_4$, which can be satisfied at the expense of the
cubic component), so that the argument in \cite[Remark 1.2]{CoMe}
works for quartics of the second type, too.

The class $4H_{\mathbb P}-2{\mathbb E}_1-2{\mathbb E}_2$ (and its
restriction onto $V_2$) is still nef and big and the map $\rho$ of
Subsection 1.1 is injective. The proceedings of Subsection 1.2
work without changes, except for the only one which is obvious:
the quadric $S\subset {\mathbb P}^4$ simply contains the two
points $o_1$, $o_2$ (and the linear system $|2H_{\mathbb
P}-{\mathbb E}_1-{\mathbb E}_2|$ is free). Word for word the same
arguments as in Subsection 1.2 reduce the factoriality of the
quartic $V$ to proving Lemma 1.1, the statement of which does not
change. In the arguments of Subsection 1.3 one has to change only
the description of the surface $W$: now it is the strict transform
of a del Pezzo surface containing the points $o_1$, $o_2$.
Repeating the arguments of Subsection 1.3 word for word, we
complete the proof of Theorem 0.1.\vspace{0.3cm}


{\bf 1.5. A generalization.} The contents and notations of the
present subsection are independent of the rest of the paper and
are not used in the sequel. Let
$$
\pi_i\colon X_i\to X_{i-1},
$$
$i=1,\dots, k$, where $k\geqslant 3$, be a sequence of blow-ups of
the points $o_{i-1}\in X_{i-1}$, where $X_0={\mathbb P}^4$, and
${\mathbb E}_i=\pi_i^{-1}(o_{i-1})\cong {\mathbb P}^3$ be the
exceptional divisors. Let $V\subset {\mathbb P}^4$ be an
irreducible hypersurface of degree 4, $V_i\subset X_i$ its strict
transform, where
$$
V_k\sim 4H_{\mathbb P}-2\sum_{i=1}^k {\mathbb E}_i
$$
is a non-singular hypersurface. Set $Q_i=V_i\cap {\mathbb E}_i$
and assume that $Q_i$ is an irreducible quadric in ${\mathbb E}_i$
for all $i=1,\dots,k$.

{\bf Theorem 1.1.} {\it Assume that the divisorial class
$$
\Delta=2H_{\mathbb P}-\sum_{i=1}^k {\mathbb E}_i
$$
is nef and big and the linear system $|\Delta|$ is free (has no
base points). Then the variety $V$ is factorial and}
$$
\mathop{\rm Pic} V_k= {\mathbb Z}H\bigoplus
\mathop{\oplus}\limits_{i=1}^k {\mathbb Z} Q_i.
$$

{\bf Proof:} repeat the arguments of Subsections 1.1 -- 1.4 word
for word.


\section{Birational rigidity}

In this section we reduce the claims of Theorems 0.2 and 0.3 to
certain facts about multiplicities of a mobile linear system
$\Sigma$ on the quartic $V$, which will be shown in \S\S 3,5 for
quartics of the first type and in \S 4 for quartics of the second
type.\vspace{0.3cm}

{\bf 2.1. Maximal singularities.} Up to, and including, Subsection
2.3, we consider quartics of the first type. Let $\Sigma$ be a
mobile linear system on $V$ (that is, $\mathop{\rm codim}
\mathop{\rm Bs} \Sigma\geqslant 2$). There is a positive integer
$n\geqslant 1$ and non-negative integers $m_1, m_2\in {\mathbb
Z}_+$, such that
$$
\Sigma\subset |nH-m_1Q_1-m_2Q_2|.
$$
(It is easy to see that $m_1\geqslant m_2$ and $n=n[\Sigma]$ in
the sense of Subsection 0.3.) Let $L_1,\dots, L_{24}$ be the lines
on $V$, passing through the point $o_1$. Repeating the arguments
of \cite[\S 1]{Pukh2025a} (which, in their turn, repeat the
arguments of \cite[\S\S 1,2]{Pukh89c}) word for word, we see that
both Theorems 0.2 and 0.3 follow from a more technical claim.

{\bf Theorem 2.1.} {\it Assume that for a general divisor
$D\in\Sigma$ the pair $(V,\frac{1}{n} D)$ is not canonical. Then
among the 25 integers
$$
m_1,\quad \mathop{\rm mult}\nolimits_{L_1}\Sigma,\quad \dots, \quad
\mathop{\rm mult}\nolimits_{L_{24}}\Sigma
$$
precisely one is strictly higher than $n$, and the others are
strictly smaller than $n$.}

(See \cite[Theorem 1.1]{Pukh2025a}.)

Here the word for word repeating comes to an end: the fact that
the vertex $o_1$ of the cone $Q_1$ is a singular point of $V_1$
changes some arguments in the proof of Theorem 2.1 compared to the
proof of \cite[Theorem 1.1]{Pukh2025a}. Let us describe in detail,
precisely which cases require additional work, and which are
completely similar to \cite{Pukh2025a}.

Let $E$ be an exceptional divisor over $V$ which is a maximal
singularity of the system $\Sigma$, that is, satisfies the
Noether-Fano inequality
$$
\mathop{\rm ord}\nolimits_E\Sigma >n\cdot a(E).
$$
The arguments of \cite[Subsection 2.1]{Pukh2025a} (that go back to
the classical paper \cite{IM}) exclude the option that the centre
of the singularity $E$ on $V$ is a non-singular point $p\neq o_1$.
Therefore, for the centre of the singularity $E$ we have two
options:

(1) this centre is the point $o_1$,

(2) this centre is an irreducible curve $B\subset V$.

In its turn, each of these options generates a few cases: for (1)
they are

(1.1) the exceptional cone $Q_1$ is a maximal singularity, that is,
$m_1>n$,

(1.2) the inequality $m_1\leqslant n$ holds, and the centre of $E$
on $V_1$ is the point $o_2\in Q_1$,

(1.3) $m_1\leqslant n$, and the centre of $E$ on $V_1$ is a point
$p\in Q_1$, $p\neq o_2$,

(1.4) $m_1\leqslant n$, and the centre of $E$ on $V_1$ is a
generator $R\ni o_2$ of the cone $Q_1$,

(1.5) $m_1\leqslant n$, and the centre of $E$ on $V_1$ is an
irreducible curve of degree $\geqslant 2$ on $Q_1$;

\noindent for (2) they are

(2.1) the centre of $E$ is one of the lines $L_1,\dots, L_{24}$,

(2.2) the centre of $E$ is a curve $B\subset V$, not containing
the point $o_1$,

(2.3) the centre of $E$ is a curve $B$ of degree $\geqslant 2$,
containing the point $o_1$.

In order to prove Theorem 2.1, we have to investigate all these 8
cases. Note that essentially this theorem claims that either (1.1)
or (2.1) takes place, and in the latter case only one line is the
centre of a maximal singularity.\vspace{0.3cm}


{\bf 2.2. Maximal curves.} Assume that the centre of the maximal
singularity $E$ is an irreducible curve $B\subset V$. Arguing word
for word as in \cite[Subsection 2.2]{Pukh2025a}, we exclude the
case (2.2) (because in that case we may just refer to \cite{IM})
and conclude that $\mathop{\rm deg} B\leqslant 3$. If the curve
$B$ is not contained in a  (two-dimensional) plane, then
$\mathop{\rm deg} B=3$ and $B$ is a rational normal curve in some
hyperplane in ${\mathbb P}^4$.

{\bf Proposition 2.1.} {\it The following inequality holds:}
$$
\mu=\mathop{\rm mult}\nolimits_B\Sigma < n.
$$

{\bf Proof.} We argue as in \cite[Subsection 2.2]{Pukh2025a}, but
with some changes coming from the singular point $o_2\in V_1$. Let
$B^1\subset V_1$ be the strict transform of $B$. If $o_2\not\in
B^1$, then for a general quadric $W$ in ${\mathbb P}^4$,
containing the curve $B$, we have $W^1\not\ni o_2$, so that the
strict transform of the surface $S=W\cap V$, that is, $S^1\subset
V_1$, is non-singular: we get the claim, similar to \cite[Lemma
2.1]{Pukh2025a} (where the strict transform of the surface $S$ is
denoted by the symbol $S^+$), so that the argument, that follows
immediately after the proof of \cite[Lemma 2.1]{Pukh2025a}, can be
repeated word for word, which gives the claim of our proposition.

Therefore we assume that $o_2\in B^1$. In that case for a general
quadric $W\supset B$ we have $\mathop{\rm Sing} S^1=\{o_2\}$.

{\bf Lemma 2.1.} {\it The point $o_2$ is a non-degenerate
quadratic singularity of the surface $S^1$, so that its strict
transform $S^2\subset V_2$ is a non-singular surface, and the
exceptional curve $S^2\cap {\mathbb E}_2$ is a non-singular conic
(a plane section of the quadric $Q_2$).}

{\bf Proof of the lemma.} The strict transforms on $X_1$ of the
quadrics in ${\mathbb P}^4$, cutting out the curve $B$, cut out
$B^1$ on $X_1$. Therefore, $W_2\cap {\mathbb E}_2$ is a general
plane in ${\mathbb E}_2\cong {\mathbb P}^3$, containing the point
$B^2\cap {\mathbb E}_2$ (when we add the upper index 2, we mean
taking the strict transform on $X_2$), so that
$$
W^2\cap Q_2=S^2\cap Q_2=(W^2\cap {\mathbb E}_2)\cap Q_2
$$
is a section of the non-singular quadric $Q_2$ by a plane of
general position in ${\mathbb E}_2$, containing the point $B^2\cap
Q_2$. Therefore, $S^2\cap Q_2$ is a non-singular conic. Q.E.D. for
the lemma.

It is easy to check that the self-intersection of the mobile part
of the linear system $\Sigma^2|_{S^2}$ is
$$
8n^2-6n\mu-4\mu^2-(2m_1-m_2)^2-(\mu-m_2)^2\geqslant 0,
$$
which implies that $\mu<n$. Q.E.D. for the proposition.

(In fact, the computations performed above in the singular case
and the computations in \cite[Subsection 2.2]{Pukh2025a} in the
non-singular case give a stronger estimate
$\mu<\frac14(-3+\sqrt{41})$, but we do not need it.)

Thus the curve $B$ is contained in some two-dimensional plane
$P\ni o_1$ in ${\mathbb P}^4$.

The plane curve $V\cap P$ is a union of irreducible components
$C_1,\dots, C_k$ of degrees $d_1,\dots, d_k$, where $k\leqslant 4$.

{\bf Proposition 2.2 (the plane section lemma).} (i) {\it If
$d_i\geqslant 2$, then
$\mathop{\rm mult}\nolimits_{C_i}\Sigma\leqslant n$.

{\rm (ii)} If $d_i=d_j=1$, where $i\neq j$, then}
$$
\mathop{\rm mult}\nolimits_{C_i}\Sigma+
\mathop{\rm mult}\nolimits_{C_j}\Sigma\leqslant 2n.
$$

{\bf Proof} will be given in \S 3.

Due to the presence of the singular point $o_2\in V_1$ we can not
(as it was done in \cite[Subsection 2.3]{Pukh2025a}) simply give a
reference to \cite[\S 6]{Pukh89c}: there are cases that are not
covered by the arguments in that paper, and they should be studied
one by one. We put this work off until \S 3.

Proposition 2.2, (i) excludes the option (2.3), and the part (ii)
provides the uniqueness of a maximal line in Theorem 2.1. To
conclude, let us complete the proof of the claim about uniqueness
in that theorem: let $L=L_i$, $1\leqslant i\leqslant 24$, be one
of the lines and $\mu=\mathop{\rm mult}\nolimits_L\Sigma$. Since
$L^1\not\ni o_2$, where $L^1\subset V_1$ is the strict transform
of that line, we conclude that the strict transform $S^1\subset
V_1$ of the section $S$ of the quartic $V$ by a general
hyperplane, containing the line $L$, does not contain the point
$o_2$, that is, $S^1\cap Q_1=S^1\cap {\mathbb E}_1$ is a
non-singular conic. Therefore, we can repeat the claim and the
proof of \cite[Lemma 2.2]{Pukh2025a} word for word: if one of the
integers $m_1$, $\mu$ is strictly higher than $n$, then the other
one is strictly smaller than $n$, which gives the full claim about
uniqueness in Theorem 2.1.\vspace{0.3cm}


{\bf 2.3. Maximal singularities over the point $o_1$.} Now let us
consider the options (1.2)--(1.5). By the work that was done in
Subsection 2.2, we may assume that the inequality $\mathop{\rm
mult}\nolimits_C\Sigma\leqslant n$ holds for every irreducible
curve $C\subset V$, and also that $m_1\leqslant n$.

{\bf Proposition 2.3.} {\it The case (1.2) does not take place:
the centre of $E$ on $V_1$ is not the point $o_2$.}

{\bf Proof.} Assume the converse. Since $m_1\leqslant n$, the pair
$(V_1,\frac{1}{n} D^1)$ is not canonical, where $D^1\in\Sigma^1$
is the strict transform of a general divisor $D^1$. Applying the
inversion of adjunction (in the same way as it was done in
\cite{Co00} or \cite{Me04}), we conclude that already $Q_2$ is a
non-canonical singularity of that pair, that is $m_2>n$. But then
$m_1\geqslant m_2>n$, too. Contradiction. Q.E.D. for the
proposition.

Since for every line $L\subset V$, such that $L\ni o_1$, we have
$L^1\not\ni o_2$, the arguments in \cite[\S 3]{Pukh2025a} work
without any changes, so that the case (1.3) does not take place.
Furthermore, the proof of \cite[Lemma 2.3]{Pukh2025a} works
without changes, so that the case (1.5) also does not take place.
However, the proof of \cite[Lemma 2.4]{Pukh2025a} does not work,
precisely because the vertex $o_2$ of the cone $Q_1$ is a
singularity of the variety $V_1$.

{\bf Proposition 2.4.} {\it The case (1.4) does not take place:
the centre of the maximal singularity $E$ on $V_1$ is not a
generator} $R\subset Q_1$.

{\bf Proof} will be given in \S 5. It requires the technique that
is used in the proof of the plane section lemma, and for that
reason makes an essential use of the work in \S 3.

This completes the proof of Theorem 2.1 for quartics of the first
type.\vspace{0.3cm}


{\bf 2.4. Quartics of the second type.} Now let $V$ be a quartic
of the second type. Here we point out the necessary changes
(mainly simplifications), that are needed in the arguments of
Subsections 2.1 -- 2.3. Again
$$
\Sigma\subset |nY-m_1Q_1-m_2Q_2|
$$
is a mobile linear system, however, the integers $m_1$, $m_2$ are
{\it not} related by the inequality $m_1\geqslant m_2$. Let
$L_{1,1},\dots, L_{1,24}$ be the lines on $V$, passing through the
point $o_1$, and $L_{2,1},\dots, L_{2,L(o_2)}$ the lines on $V$,
passing through the point $o_2$. Theorem 2.1 should be replaced by
a claim about $26+L(o_2)$ multiplicities.

{\bf Theorem 2.2.} {\it Assume that for a general divisor
$D\in\Sigma$ the pair $(V,\frac{1}{n} D)$ is not canonical. In
that case among $26+L(o_2)$ multiplicities
$$
m_1,\quad \mathop{\rm mult}\nolimits_{L_{1,i}}\Sigma,
1\leqslant i\leqslant 24,\quad m_2,\quad
\mathop{\rm mult}\nolimits_{L_{2,i}}\Sigma, 1\leqslant i\leqslant L(o_2),
$$
precisely one is strictly higher than $n$, and the others are
strictly smaller than} $n$.

Fix a maximal singularity $E$ of the system $\Sigma$. Assume that
its centre on $V$ is an irreducible curve $B\subset V$. Again
$\mathop{\rm deg} B\leqslant 3$ and if $B$ is not contained in a
two-dimensional plane, then $B$ is a rational normal curve in some
hyperplane in ${\mathbb P}^4$. If $B$ does not contain the points
$o_1$, $o_2$, then we exclude this case, just making a reference
to \cite{IM}; if $B$ contains only one of them, then we exclude
this case, referring to \cite{Pukh89c,Pukh2025a}. It remains to
check that the claim of Proposition 2.1 is true if $B$ contains
both points $o_1$, $o_2$. Arguing as in the proof of Proposition
2.1, consider the surface $S\supset B$, cut out on $V$ by a
quadric, such that its strict transform $S^2\subset V_2$ is
non-singular and the restriction $\varphi_1\circ\varphi_2\colon
S^2\to S$ is the blow-up of two non-degenerate quadratic
singularities $o_1$, $o_2$. The restriction $\Sigma^2|_{S^2}$ has
a unique fixed component, the curve $B^2$ of multiplicity $\mu$.
The self-intersection of the mobile part of this system is
$$
8n^2-6n\mu-3\mu^2-m_1^2-m_2^2-(\mu-m_1)^2-(\mu-m_2)^2\geqslant 0,
$$
which excludes $B$ as the centre of the maximal singularity $E$.

The plane section lemma (Proposition 2.2) remains true.

If the plane $P$ contains only one of the two singular points,
then it is sufficient to refer to the proof of the plane section
lemma in \cite[\S 6]{Pukh89c}. If $P\ni o_1,o_2$, then this proof
does not work for formal reasons and the claim of the plane
section lemma is shown in \S 4.

Now in order to prove the claim about uniqueness in Theorem 2.2,
it remains to show the following three facts.

{\bf Proposition 2.4.} {\it If one of the two integers $m_1, m_2$
is strictly higher than $n$, then the other one is strictly
smaller than} $n$.

{\bf Proposition 2.5.} {\it Let $\{i,j\}=\{1,2\}$. If one of the
two integers $m_i$ and $\mathop{\rm
mult}\nolimits_{L_{j,l}}\Sigma$ is strictly higher than $n$, then
the other one is strictly smaller than} $n$.

{\bf Proposition 2.6.} {\it Let $L_{1,i}$ and $L_{2,j}$ be two
lines on $V$, where $L_{1,i}\cap L_{2,j}=\emptyset$. Then if one
of the two integers
$$
\mathop{\rm mult}\nolimits_{L_{1,i}}\Sigma,\quad
\mathop{\rm mult}\nolimits_{L_{2,j}}\Sigma
$$
is strictly higher than $n$, then the other one is strictly
smaller than} $n$.

{\bf Proof of Proposition 2.4.} Consider the section $S$ of the
hypersurface $V$ by a general hyperplane, containing the line
$[o_1,o_2]$. Since $[o_1,o_2]\not\subset V$, the restriction
$\Sigma|_S$ has no fixed components, so that $\Sigma^2|_{S^2}$
does not have them, too. The self-intersection of the latter
linear system is
$$
4n^2-2m_1^2-2m_2^2\geqslant 0,
$$
which completes the proof of the proposition.

{\bf Proof of Proposition 2.5} we put off until \S 4, since the
point $o_i$ and the line $L_{j,l}$ are contained in some plane,
which requires the proof of the plane section lemma.

{\bf Proof of Proposition 2.6.} Set
$$
\mu_1=\mathop{\rm mult}\nolimits_{L_{1,i}}\Sigma,\quad
\mu_2=\mathop{\rm mult}\nolimits_{L_{2,j}}\Sigma.
$$
Consider the section $S$ of the hypersurface $V$ by a general
hyperplane, containing the line $L_{1,i}$. The restriction
$\Sigma|_S$ has the unique fixed component $L_{1,i}$ of
multiplicity $\mu_1$ and a base point $S\cap L_{2,j}$ of
multiplicity $\mu_2$. Computing the self-intersection of the
mobile part of the restriction $\Sigma^1|_{S^1}$ (see \cite[Lemma
2.2]{Pukh2025a}), we obtain the inequality
$$
4n^2-m_1^2-\mu_1^2-2n\mu_1-(m_1-\mu_1)^2\geqslant \mu_2^2,
$$
which completes the proof of the proposition.

Now in order to prove Theorem 2.2 it remains to show that if all
$26+L(o_2)$ multiplicities do not exceed $n$, then the assumption
that there exists a maximal singularity $E$ leads to a
contradiction. By what has already been done, the centre of $E$ on
$V$ is either the point $o_1$, and this option is excluded in
\cite[\S 3]{Pukh2025a}, or the point $o_2$, and that option is
excluded by means of the inversion of adjunction \cite{Co00,Me04}.
Proof of Theorem 2.2 is complete. Therefore, Theorems 0.2 and 0.3
are shown for quartics of the second type.


\section{The plane section lemma. I}

In this section we prove the plane section lemma (Proposition 2.2)
for quartics of the first type. Apart from that, we will show a
certain claim about the multiplicity of the strict transform of
the system $\Sigma$ on $V_1$ along a generator of the cone $Q_1$:
the technique of the proof of that claim (which will be needed in
\S 5 in the proof of Proposition 2.4) is close to the proof of the
plane section lemma.\vspace{0.3cm}

{\bf 3.1. Statement of the problem.} Recall that we consider
quartics of the first type, so that $o_2$ is the vertex of the
cone $Q_1$. Let $P\ni o_1$ be some two-dimensional plane,
$P^1\subset X_1$ its strict transform and $R=P^1\cap{\mathbb E}_1$
the exceptional line of the blow up of the point $o_1$ on $P$. If
$R$ is not a generator of the cone $Q_1$, then $R\not\subset Q_1$.
Assume first that $o_2\not\in R$, that is, $o_2\not\in P^1$. In
that case the proof of the plane section lemma in \cite[\S
6]{Pukh89c} works without changes so that the plane section lemma
is shown (see also the comments in \cite[Proposition
2.1]{Pukh2025a}).

Assume that $o_2\in P^1$. Set
$$
\mu=\mathop{\rm mult}\nolimits_{o_1}(P\circ V),
$$
where $(P\circ V)$ is the plane quartic cut out by the
hypersurface $V$ on the plane $P$, it can be reducible and contain
multiple components. Its affine part is given by the equation
$f(z_*)|_P=0$. It is easy to see that $R\subset(P^1\circ V^1)$ if
and only if $\mu\geqslant 3$, and $R$ comes into the effective
1-cycle $(P^1\circ V^1)$ with multiplicity $\mu - 2$.

In order to reduce the number of cases which we have to consider,
assume in addition that the quartic $V$ satisfies the following
condition of general position: in the notations of Subsection 0.2
the polynomial $q_3(z_*)$ does not vanish identically on all
generators of the cone $Q_1=\{q_2=0\}$ (quartics, for which this
condition is violated, form a subset of codimension 2 in the
quasi-projective variety of quartics of the first type, and
moreover, for every plane $P\ni o_1$ the inequality $\mu\leqslant
3$ holds, that is, the effective cycle $(P\circ V)$ is neither
$a_iL_i+a_jL_j$, where $a_i+a_j=4$ and $i\neq j$, nor $4L_i$ (here
$i,j\in\{1,\dots,24\}$). Therefore, $\mu\in\{2,3\}$.

Let us list all options for the cycle $(P\circ V)$. Each option
has its number in the format (a.b), where $a=\mu$. We consider the
planes, for which $(P\circ V)$ is not an irreducible curve of
degree 4, such curves can not be the centre of a maximal
singularity, see Subsection 2.2.

(2.1) $(P\circ V)=C_1+C_2$, where $C_1,C_2$ are irreducible
conics, tangent at the point $o_1$, and moreover, $o_2\in
C^1_1,C^1_2$.

(2.2) $(P\circ V)=2C$, where $C$ is an irreducible conic, such
that $o_2\in C^1$.

(This completes the list of options for $\mu=2$. It may seem that
there are many more options, but this is not the case because of
the assumption that $[o_1,o_2]\not\subset V$ and the condition
$\mathop{\rm mult}_{o_2}(P^1\circ V_1)=2$, taking into account
that $R$ is not a generator of the cycle $(P^1\circ V_1)$.)

(3.1) $(P\circ V)=C_1+C_2$, where $C_1$ is an irreducible cubic
with the double point $o_1$, $C_2\ni o_1$ is a line, and moreover,
$o_2\in C^1_1$ and $o_2\not\in C^1_2$. The cycle $(P^1\circ V_1)$
is of the form $C_1^1+C_2^1+R$.

(3.2) $(P\circ V)=C_1+C_2+C_3$, where $C_1$ is an irreducible
conic, $C_2$ and $C_3$ are lines, and moreover, $o_2\in C^1_1$ and
$o_2\not\in C^1_2$, $o_2\not\in C^1_3$. The cycle $(P^1\circ V_1)$
is of the form $C^1_1+C^1_2+C^1_3+R$.

(3.3) $(P\circ V)=C_1+2C_2$, where $C_1$ is an irreducible conic
and $C_2$ is a line, $o_2\in C^1_1$.

The list above is complete. It is easy to see that for a quartic
of the first type of general position, of the five cases listed
above only (3.1) takes place for some plane $P$ (precisely because
we do not need to consider planes, intersecting $V$ by an {\it
irreducible} curve of degree 4). We will prove the plane section
lemma for the cases (2.1), (3.1) and (3.2), in order to show that
the explicit geometric version of the technique of the proof of
that lemma that was used in \cite{Pukh89c} works for all cases.
The quartics, for which (2.2) or (3.3) take place, form subsets of
high codimension in the set of quartics of the first
type.\vspace{0.3cm}


{\bf 3.2. The case of a cubic curve and a line.} Assume that the
case (3.1) takes place. Let $S$ be the section of $V$ by a general
hyperplane, containing the plane $P$.

{\bf Lemma 3.1.} {\it The point $o_1$ is the only singularity of
the surface $S$ and the point $o_2$ is the only singularity of its
strict transform $S^1$, and the strict transform $S^2$ is a
non-singular surface.}

{\bf Proof.} Consider the pencil  $|H-P|$ of sections of the
hypersurface $V$ by hyperplanes in ${\mathbb P}^4$, containing the
plane $P$. The intersection $(S_1\circ S_2)$ of two general
distinct elements of that pencil is, obviously, the 1-cycle
$(P\circ V)$. Since $C_1\neq C_2$ are components of multiplicity
one of the latter intersection, the possible singularities of the
surfaces $S_1$, $S_2$ are only the points of intersection $C_1\cap
C_2$. For any point $p\in C_1\cap C_2$, which is not the point
$o_1$, the surface $S$ has a singularity at the point $p$ if and
only if $S$ is cut out by the tangent hyperplane $T_pV$.
Therefore, for a general surface $S\in|H-P|$ we have $\mathop{\rm
Sing} S=\{o_1\}$.

Hyperplanes in ${\mathbb P}^4$, containing the plane $P$,
correspond to planes in ${\mathbb E}_1\cong {\mathbb P}^3$,
containing the line $R$. Therefore we get:
$$
S^1\cap {\mathbb E}_1=S^1\cap E_1=R_1\cup R_2,
$$
where $R_1=R$ and $R_2$ is the other generator of the cone $Q_1$,
varying together with $S$. This is true at the scheme-theoretic
level, too,
$$
(S^1\circ {\mathbb E}_1)=R_1+R_2,
$$
so that $o_2=R_1\cap R_2$ is the only singularity of the surface
$S^1$.

Finally, let $P^2\subset X_2$ be the strict transform of the plane
$P$. Since
$$
\mathop{\rm mult}\nolimits_{o_2} (P^1\circ V_1)=2
$$
(the strict transform of the cubic $C_1$ is non-singular and
$C_2^1\not\ni o_2$), we get
$$
(P^2\circ V_2)=C_1^2+C_2^2+R^2,
$$
that is, this scheme-theoretic intersection does not contain the
exceptional line $P^2\cap {\mathbb E}_2$, which is the exceptional
curve of the blow-up $\pi_2|_{P^2}\colon P^2\to P^1$ of the point
$o_2$ on $P^1$. Thus this line $P^2\cap {\mathbb E}_2$ is not
contained in the quadric $Q_2$, so that the surfaces $S^2$ for
$S\in |H-P|$ cut out on $Q_2$ the pencil of conics which is cut
out on $Q_2$ by the pencil of planes, containing the line $P^2\cap
{\mathbb E}_2$. Therefore $S^2\cap Q_2$ is a non-singular conic,
so that the surface $S^2$ is non-singular. Q.E.D. for the lemma.

Let $R_3=S^2\cap Q_2$ be the exceptional conic of the blow-up of
the quadratic singularity $o_2\in S^1$. We look at the
multiplicities of the linear system of curves
$\Sigma_S=\Sigma|_S$, pulled back on $S^2$, along the curves
$C_1^2$, $C_2^2=C_2^1$, $R_1^2$, $R_2^2$, $R_3$. Intersection
numbers give the following multiplication table for these five
curves:
\begin{equation}\label{03.11.2025.1}
\left(
\begin{array}{ccccc}
-2 & 1  & 1  & 0  & 1\\
1  & -2 & 1  & 0  & 0\\
1  & 1  & -2 & 0  & 1\\
0  & 0  & 0  & -2 & 1\\
1  & 0  & 1  & 1  & -2
\end{array}\right)
\end{equation}
(the columns and rows correspond to the order, in which these
curves are listed above). The canonical class of the surfaces $S$,
$S^1$, $S^2$ is equal to zero. The pull back of the hyperplane
section $H_S$ on $S^2$ (denoted by the same symbol) is
$$
H_S=C_1^2+C_2^2+2R_1^2+R_2^2+2R_3,
$$
the class of the restriction $E_{1,S}={\mathbb
E}_1|_{S^1}=Q_1|_{S^1}$, pulled back on $S^2$, is
$$
E_{1,S}=R_1^2+R_2^2+R_3,
$$
and $E_{2,S}={\mathbb E}_2|_{S^2}=Q_2|_{S^2}=R_3$. The class of
the mobile part of the restriction of the strict transform
$\Sigma^2$ of the linear system $\Sigma$ onto the surface $S^2$ is
$$
\xi=nH_S-m_1E_{1,S}-m_2E_{2,S}-bR_1^2-b_1C_1^2-b_2C_2^2=
$$
$$
=(n-b_1)C_1^2+(n-b_2)C_2^2+(2n-m_1-b)R_1^2+(n-m_1)R_2^2+
(2n-m_1-m_2)R_3.
$$
Clearly, this class has a non-negative intersection with each of
the five curves. Now consider the ${\mathbb Q}$-vector space
$$
\Xi={\mathbb Q}C_1^2\oplus {\mathbb Q}C_2^2\oplus {\mathbb Q}R_1^2
\oplus {\mathbb Q}R_2^2\oplus {\mathbb Q}R_3
$$
and define the new classes $C_1^+$, $C_2^+$, $R_1^+$ by the formulas
$$
\begin{array}{ccl}
C_1^2 & = & C_1^2 +\frac13 R_2^2 +\frac23 R_3,\\
C_2^+ & = & C_2^2,\\
R_1^+ & = & R_1^2 +\frac13 R_2^2 +\frac23 R_3.
\end{array}
$$
(Obviously, the intersection form on $\Xi$ is given by the
non-degenerate matrix (\ref{03.11.2025.1}).) Set
$$
\Xi^+={\mathbb Q}C_1^+\oplus {\mathbb Q}C_2^+\oplus
{\mathbb Q}R_1^+.
$$
The subspace $\Xi^+$ is orthogonal to the classes $R_2^2$ and
$R_3$. The intersection form on this space is given by the matrix
\begin{equation}\label{03.11.2025.2}
\left(
\begin{array}{ccc}
-\frac43 & 1  & \frac53\\
1        & -2 & 1      \\
\frac53  & 1  & -\frac43
\end{array}\right).
\end{equation}
Set $\xi^+=\xi_1C_1^++\xi_2C_2^++\xi_3R_1^+$, where $\xi_1=n-b_1$,
$\xi_2=n-b_2$, $\xi_3=2n-m_1-b$.

{\bf Lemma 3.2.} {\it The intersections $(\xi^+\cdot C_1^+)$,
$(\xi^+\cdot C_2^+)$, $(\xi^+\cdot R_1^+)$ are non-negative.}

{\bf Proof.} Since $C_1^+$, $C_2^+$, $R_1^+$ are orthogonal to
$R_2^2$ and $R_3$, we get the equalities $(\xi^+\cdot
C_1^+)=(\xi\cdot C_1^+)$, $(\xi^+\cdot C_2^+)=(\xi\cdot C_2^+)$
and $(\xi^+\cdot R_1^+)=(\xi\cdot R_1^+)$. However, the right hand
sides of these equalities are non-negative, because the classes
$C_1^+$, $C_2^+$, $R_1^+$ are non-negative linear combinations of
the classes of the standard basis of the space $\Xi$, and $\xi$
has a non-negative intersection with these classes. Q.E.D. for the
lemma.

Let $\Theta=\|\theta_{ij}\|_{1\leqslant i,j\leqslant 3}$ be the
inverse matrix for (\ref{03.11.2025.2}). We have the obvious
equalities
$$
\xi_i=\theta_{i1} (\xi^+\cdot C_1^+)+\theta_{i2} (\xi^+\cdot C_2^+)+
\theta_{i3} (\xi^+\cdot R_1^+),
$$
where $i=1,2,3$. Lemma 3.2 implies that if the $i$-the row (or the
$i$-th column, the matrix $\Theta$ is symmetric) consists of
non-negative numbers, then $\xi_i\geqslant 0$. It is easy to check
that
$$
\Theta=\frac{1}{24}\left(
\begin{array}{ccc}
5  & 9  & 13 \\
9  & -3 & 9  \\
13 & 9  & 5
\end{array}\right).
$$
This immediately implies the inequality $\xi_1\geqslant 0$, which
is equivalent to the estimate
$$
\mathop{\rm mult}\nolimits_{C_1}\Sigma\leqslant n.
$$
This is all that Proposition 2.2 claims in the case (3.1), so that
the plane section lemma is shown for that case.

Note that $\theta_{22}=-\frac{3}{24}$, and this agrees with the
fact that the line $C_2$ can be a maximal curve of the system
$\Sigma$.

In order to prove Proposition 2.2 it was sufficient to consider
the orthogonal complement to the three classes $R_1^2$, $R_2^2$,
$R_3$. We did not do that intentionally: from the explicit form of
the matrix $\Theta$ we conclude that $\xi_3\geqslant 0$, so that
the following claim is true.

{\bf Proposition 3.1.} {\it The following inequality holds}
$$
m_1+\mathop{\rm mult}\nolimits_{R}\Sigma^1\leqslant 2n.
$$

Geometrically it means that the exceptional divisor of the blow-up
of the infinitely near line $R$ on $V_1$ is not a maximal
singularity of the system $\Sigma$. This fact will be needed in \S
5 in the proof of Proposition 2.4.

{\bf Remark 3.1.} The argument above, when we restrict the linear
system $\Sigma$ onto a general surface $S\in |H-P|$ and study the
intersection form on the vector space, the basis of which is
formed by the numerical classes of the components of the plane
section and the exceptional curves of the resolution of
singularities of the surface $S$, is the explicit geometric form
of the argument in \cite[\S 6]{Pukh89c}, which was carried out in
terms of numerical classes of a certain blow-up of the quartic. In
\cite[\S 6]{Pukh89c} it was done in that way in order to avoid
considering a large number of possible cases. In this paper we
restrict the amount of cases by the conditions of general position
and so it is easier to prove the plane section lemma for each
case, in which the point $o_2$ is present, by explicit geometric
arguments, resolving singularities of the surface $S$ and studying
the intersection form for a particular set of curves. However, the
logic of the arguments is the same as in \cite[\S 6]{Pukh89c}; in
particular, the replacement of the classes $C_1^2$, $C_2^2$,
$R_1^2$ by the classes $C_1^+$, $C_2^+$, $R_1^+$, orthogonal to
the curves $R_2^2$, $R_3$, corresponds to the ``reductions'' in
\cite[\S 6]{Pukh89c}.

For the remaining cases (2.1) and (3.2) we only point out the
necessary changes in the arguments for the case (3.1) considered
above.\vspace{0.3cm}


{\bf 3.3. The case of a conic and two distinct lines.} Assume that
the case (3.2) takes place. The claim of Lemma 3.1 remains true
(the proof can be repeated word for word). On the non-singular
surface $S^2$ in this case we have 6 curves
$$
C_1^2,\quad C_2^2=C_2^1,\quad C_3^2=C_3^1,\quad R_1^2,\quad R_2^2,
\quad R_3,
$$
the intersection indices of which form the matrix
$$
\left(\begin{array}{cccccc}
-2 & 1  & 1  & 0  & 0  & 1  \\
1  & -2 & 0  & 1  & 0  & 0  \\
1  & 0  & -2 & 1  & 0  & 0  \\
0  & 1  & 1  & -2 & 0  & 1  \\
0  & 0  & 0  & 0  & -2 & 1  \\
1  & 0  & 0  & 1  & 1  & -2
\end{array}\right).
$$
The pull back of the hyperplane section $H_S$ on $S^2$ is
$$
H_S=C_1^2+C_2^2+C_3^2+2R_1^2+R_2^2+2R_3,
$$
the classes $E_{1,S}$ and $E_{2,S}$ are given by the same formulas
as in the case (3.1), so that the class of the mobile part of the
restriction of the strict transform of $\Sigma^2$ onto the surface
$S^2$ is
$$
\xi=\sum_{i=1}^3(n-b_i)C_i^2+(2n-m_1-b)R_1^2+(n-m_1)R_2^2+
(2n-m_1-m_2)R_3,
$$
where $b_i=\mathop{\rm mult}\nolimits_{C_i}\Sigma$ and
$b=\mathop{\rm mult}\nolimits_{R}\Sigma^1$. By construction, $\xi$
intersects non-negatively each of the 6 curves. Consider the
vector space
$$
\Xi={\mathbb Q}C_1^2\oplus {\mathbb Q}C_2^2\oplus
{\mathbb Q}C_3^2\oplus
{\mathbb Q}R_1^2 \oplus {\mathbb Q}R_2^2\oplus {\mathbb Q}R_3
$$
and define the new classes $C_1^+$, $C_2^+$, $C_3^+$, $R_1^+$,
orthogonal to $R_2^2$ and $R_3$, by the formulas
$$
\begin{array}{ccl}
C_1^2 & = & C_1^2 +\frac13 R_2^2 +\frac23 R_3,\\
C_2^+ & = & C_2^2,\\
C_3^+ & = & C_3^2,\\
R_1^+ & = & R_1^2 +\frac13 R_2^2 +\frac23 R_3,
\end{array}
$$
so that the intersection form for the basis classes of the vector
space
$$
\Xi^+={\mathbb Q}C_1^+\oplus {\mathbb Q}C_2^+\oplus {\mathbb
Q}C_3^+\oplus {\mathbb Q}R_1^+
$$
is given by the matrix
\begin{equation}\label{04.11.2025.1}
\left(
\begin{array}{cccc}
-\frac43 & 1  & 1  & \frac23\\
1        & -2 & 0  & 1      \\
1        & 0  & -2 & 1      \\
\frac23  & 1  & 1  & -\frac43
\end{array}\right).
\end{equation}
Set $\xi^+=\sum_{i=1}^3 \xi_i C_i^++\xi_4 R_1^+$, where for
$i\in\{1,2,3\}$ we have $\xi_i=n-b_i$ and $\xi_4=2n-m_1-b$.
Repeating the arguments in the proof of Lemma 3.2 word for word,
we conclude that $(\xi^+\cdot C_i^+)\geqslant 0$ for
$i\in\{1,2,3\}$ and $(\xi^+\cdot R_1^+)\geqslant 0$. The inverse
matrix for (\ref{04.11.2025.1}) is
$$
\Theta=\|\theta_{ij}\|_{1\leqslant i,j\leqslant 4}=
\frac{1}{8}\left(
\begin{array}{cccc}
1  & 3  & 3  & 5 \\
3  & -1 & 3  & 3 \\
3  &  3 & -1 & 3 \\
5  &  3 & 3  & 1
\end{array}\right).
$$
As in the case (3.1), the non-negativity of the numbers in the
first row of the matrix $\Theta$ gives the first claim of
Proposition 2.2. The equality
$$
\xi_2+\xi_3=\sum_{i=1}^3 (\theta_{2i}+\theta_{3i})(\xi^+\cdot C_i^+)+
(\theta_{24}+\theta_{34})(\xi^+\cdot R_1^+)
$$
and the fact that the sum of the second and third rows in the
matrix $\Theta$ consists of positive numbers imply that
$\xi_2+\xi_3\geqslant 0$, which is equivalent to the second claim
of Proposition 2.2. Proof of the plane section lemma in the case
(3.2) is complete.

The non-negativity (in fact, positivity, but the non-negativity is
sufficient) of the last row of the matrix $\Theta$ gives the claim
of Proposition 3.1 for the case (3.2).\vspace{0.3cm}


{\bf 3.4. The case of two conics.} Assume that the case (2.1)
takes place. The surface $S^2$ is again non-singular, and the
multiplication table for the set of 5 curves $C_1^2, C_2^2, R_1^2,
R_2^2, R_3$ is given by the matrix
$$
\left(\begin{array}{ccccc}
-2 & 2  & 0  & 0  &  1  \\
2  & -2 & 0  & 0  &  1  \\
0  & 0  & -2 & 0  & 1  \\
0  & 0  & 0  & -2 & 1  \\
1  & 1  & 1  & 1  & -2
\end{array}\right).
$$
Since in this case $\mu=2$, the exceptional line $R$ is not
contained in the intersection $P^1\cap V_1$, so that for the pull
back on $S^2$ of the hyperplane section of the surface $S$ we have
$$
H_S=C_1^2+C_2^2+R_1^2+R_2^2+2R_3
$$
(where, in contrast to the previous cases, $R_1\neq R$: the line
$R\subset {\mathbb E}_1$ intersects the cone $Q_1$ at the only
point, the vertex $o_2$), besides,
$$
E_{1,S}=R_1^2+R_2^2+R_3\quad\mbox{and}\quad E_{2,S}=R_3.
$$
The class of restriction of the strict transform of $\Sigma^2$
onto $S^2$ is
$$
\xi=\sum_{i=1}^2(n-b_i)C_i^2+(n-m_1)R_1^2+(n-m_1)R_2^2+(2n-m_1-m_2)R_3,
$$
where $b_i$ have the same meaning as in the previous cases. The
classes
$$
C_1^+,C_2^+\in\Xi={\mathbb Q}C_1^2\oplus {\mathbb Q}C_2^2\oplus
{\mathbb Q}R_1^2 \oplus {\mathbb Q}R_2^2\oplus {\mathbb Q}R_3
$$
are defined by the equalities
$$
C_i^+=C_i^2+\frac12 R_1^2+\frac12R_2^2+R_3,
$$
where $i=1,2$, and the intersection form on $\Xi^+={\mathbb
Q}C_1^+\oplus {\mathbb Q}C_2^+$ is given by the matrix
$$
\left(
\begin{array}{cc}
-1  & 3  \\
3  & -1
\end{array}\right),
$$
the inverse of which is
$$
\Theta=\frac{1}{8}\left(
\begin{array}{cc}
1  & 3 \\
3  & 1
\end{array}\right).
$$
Since all coefficients of that matrix are positive, the plane
section lemma is shown for the case (2.1). Q.E.D. for Proposition
2.2.\vspace{0.3cm}


{\bf 3.5. The case of a quartic with a triple point.} Assume that
$P\cap V$ is an irreducible curve $C$ of degree 4, and moreover,
$\mathop{\rm mult}\nolimits_{o_1} C=3$. Although this case is not
needed for the proof of the plane section lemma (we know that the
curve $C$ can not be maximal), it is needed in order to complete
the proof of Proposition 3.1: obviously,
$$
(P^1\circ V_1)=C^1+R,
$$
where $R$ is a generator of the cone $Q_1$, here $o_2\in C^1$ is a
smooth point (the curve $C^1$ is non-singular), so that $C^1$
meets $R$ in three points, one of which is $o_2$. For a general
surface $S\in |H-P|$ its strict transform $S^2$ is non-singular,
the multiplication table for the set of curves $C^2$, $R_1^2$,
$R_2^2$, $R_3$ (where $R=R_1$) is given by the matrix
$$
\left(\begin{array}{cccc}
-2  & 2  & 0  & 1 \\
2   & -2 & 0  & 1 \\
0   &  0 & -2 & 1 \\
1   &  1 & 1  & -2
\end{array}\right).
$$
We have the following equalities:
$$
\begin{array}{ccl}
H_S     & = & C^2 +2R_1^2+R_2^2+2R_3,\\
E_{1,S} & = & R_1^2+R_2^2+R_3,\\
E_{2,S} & = & R_3.
\end{array}
$$
The classes
$$
C^+,R_1^+\in\Xi={\mathbb Q} C^2\oplus {\mathbb Q} R_1^2\oplus
{\mathbb Q} R_2^2\oplus {\mathbb Q} R_3
$$
are given by the equalities
$$
\begin{array}{ccl}
C^+   & = & C^2+\frac13 R_2^2 +\frac23 R_3,\\
R_1^+ & = & R_1^2+\frac13 R_2^2 +\frac23 R_3,
\end{array}
$$
the intersection form on $\Xi^+={\mathbb Q} C^+\oplus {\mathbb Q}
R_1^+$ is given by the matrix
$$
\frac13\left(
\begin{array}{cc}
-4 & 8 \\
8  & -4
\end{array}\right),
$$
the inverse for which is
$$
\Theta=\frac14\left(
\begin{array}{cc}
1 & 2 \\
2  & 1
\end{array}\right).
$$
The non-negativity of coefficients of that matrix implies
Proposition 3.1 for the case under consideration.


\section{The plane section lemma. II}

In this section we prove the plane section lemma (Proposition 2.2)
for quartics of the second type. After that we prove Proposition
0.1: for a quartics of the second type of general position there
are 24 distinct lines through the point $o_2$, lying on the
quartic.\vspace{0.3cm}

{\bf 4.1. Statement of the problem.} We consider quartics of the
second type with two singular points: $o_1$ of rank 3 (which is
resolved by one blow-up) and $o_2$ of rank 4. If the plane $P$
contains only one of these two points, then the proof of the plane
section lemma given in \cite[\S 6]{Pukh89c}, works word for word.
For this reason, assume that $P\ni o_1,o_2$.

We know that the line $[o_1,o_2]\not\subset V$. Assume in addition
that for every plane $P\supset [o_1,o_2]$ the effective 1-cycle
$(P\circ V)$ does not contain multiple components (quartics of the
second type, for which this condition is not satisfied, form a
subset of positive codimension in the set of quartics of the
second type). Since $(P\circ V)=P\cap V$ is singular at the points
$o_1$, $o_2$, does not contain the line $[o_1,o_2]$ and is
reducible (as we know, an irreducible curve of degree $\geqslant
4$ can not be maximal for $\Sigma$), the following cases are
possible:

(1) $P\cap V$ is the union of four distinct lines $C_1$, $C_2$,
$C_3$, $C_4$, where $o_1\in C_1, C_2$ and $o_2\in C_3, C_4$;

(2.1) $(P\circ V)=C_1+C_2+C_3$, where $C_1$ is an irreducible
conic, $C_1\ni o_1, o_2$, $C_2$ and $C_3$ are lines, containing,
respectively, $o_1$ and $o_2$;

(2.2) $(P\circ V)=C_1+C_2$, where $C_1$, $C_2$ are irreducible
conics, passing through both points $o_1$, $o_2$;

(3) $(P\circ V)=C_1+C_2$, where $C_1$ is an irreducible cubic with
a double point and $C_2$ is aline; in principle we should consider
separately the cases (3.1), when $o_1\in C_1$ is a double point
and $o_2\in C_1\cap C_2$, and (3.2), when $o_2\in C_1$ is a double
point and $o_1\in C_1\cap C_2$, however, the arguments are
identical in both cases and we will consider the case (3.1) only.

In the next subsection we prove Proposition 2.2, using the
technique that was applied in \S 3. We will just point out the
necessary changes in the arguments and computations of the
previous section.\vspace{0.3cm}


{\bf 4.2. Proof of the plane section lemma.} In the notations of
Subsection 3.2, let $S\in |H-P|$ be a general surface. Since the
cycle $(P\circ V)$ has no multiple components, we can argue as in
the proof of Lemma 3.1 and conclude that
$$
\mathop{\rm Sing} S=\{o_1,o_2\},
$$
and the strict transform $S^2\subset V_2$ is non-singular.
Obviously, $S^2\cap {\mathbb E}_1$ and $S^2\cap {\mathbb E}_2$ are
non-singular conics, we denote them by the symbols $R_1$ and
$R_2$.

{\bf The case of four lines.} Assume that the case (1) takes
place. The multiplication table for the six curves
$$
R_1,\quad C_1^2,\quad C_2^2,\quad C_3^2,\quad C_4^2,\quad R_2
$$
(this order is more convenient) is given by the (6$\times$6)-matrix
$$
\left(\begin{array}{cccccc}
-2 & 1  &  1 & 0  & 0  &  0  \\
1  & -2 &  0 & 1  & 1  &  0  \\
1  & 0  & -2 & 1  & 1  &  0  \\
0  & 1  & 1  & -2 & 0  &  1  \\
0  & 1  & 1  & 0  & -2 &  1  \\
0  & 0  & 0  & 1  & 1  &  -2
\end{array}\right).
$$
Introducing the classes
$$
C_1^+=C_1^2+\frac12 R_1,\quad C_2^+=C_2^2+\frac12 R_1,
$$
$$
C_3^+=C_3^2+\frac12 R_2,\quad C_4^+=C_4^2+\frac12 R_2,
$$
we obtain for $C_1^+$, $C_2^+$, $C_3^+$, $C_4^+$ the
multiplication table
$$
\frac12\left(
\begin{array}{cccc}
-3 & 1  & 2  & 2  \\
1  & -3 & 2  & 2  \\
2  & 2  & -3 & 1  \\
2  & 2  & 1  & -3
\end{array}\right).
$$
The inverse matrix is
$$
\Theta=\frac{1}{12}\left(
\begin{array}{cccc}
-1 & 5  & 4  & 4  \\
5  & -1 & 4  & 4  \\
4  & 4  & -1 & 5  \\
4  & 4  & 5  & -1
\end{array}\right).
$$
Since the sum of any two distinct rows of the matrix $\Theta$ is a
vector, all coordinates of which are positive, we obtain the claim
(ii) of Proposition 2.2 (the claim (i) in this case does not
apply).

{\bf The case of a conic and two lines.} Assume that the case
(2.1) takes place. Here the multiplication table for $C_1^+$,
$C_2^+$, $C_3^+$ is given by the matrix
$$
\frac12\left(
\begin{array}{ccc}
-2 & 3  & 3  \\
3  & -3 & 2  \\
3  & 2  & -3
\end{array}\right),
$$
the inverse of which is
$$
\Theta=\frac{1}{40}\left(
\begin{array}{ccc}
5  & 15 & 15 \\
15 & -3 & 13 \\
15 & 13 & -3
\end{array}\right),
$$
which implies both claims of Proposition 2.2 in the case under
consideration.

{\bf The case of two conics.} Assume that the case (2.2) takes
place. Here the multiplication table for
$$
C_1^+=C_1^2+\frac12 R_1+\frac12 R_2,\quad C_2^+=C_2^2+\frac12
R_1+\frac12 R_2
$$
is given by the matrix
$$
\left(
\begin{array}{cc}
-1 & 3 \\
3  & -1
\end{array}\right),
$$
the inverse of which is
$$
\Theta=\frac{1}{8}\left(
\begin{array}{cc}
1  & 3 \\
3  & 1
\end{array}\right),
$$
which proves the plane section lemma in this case.

{\bf The case of a cubic curve and a line.} Assume that the case
(3.1) takes place. Here the multiplication table for
$C_1^+=C_1^2+R_1+\frac12 R_2$, $C_2^+=C_2^2+\frac12 R_2$, is given
by the matrix
$$
\frac12\left(
\begin{array}{cc}
1 & 5 \\
5  & -3
\end{array}\right),
$$
the inverse of which is
$$
\Theta=\frac{1}{14}\left(
\begin{array}{cc}
3  & 5 \\
5  & -1
\end{array}\right).
$$
This completes the proof of the plane section lemma for quartics
of the second type.\vspace{0.3cm}


{\bf 4.3. Proof of Proposition 2.5.} Let us show that one of the
singular points $o_i$, $i\in \{1,2\}$, and a line $L$, passing
through the other singular point, can not maximal simultaneously.
We will show the inequality
$$
m_i+\mathop{\rm mult}\nolimits_L \Sigma\leqslant 2n,
$$
which immediately implies Proposition 2.5. Let $P=\langle o_i,
L\rangle$ be the plane, containing this point and this line. Proof
follows the scheme of the proof of the plane section lemma in the
cases (1), (2.1) and (3), see Subsection 4.1, we only need a
different orthogonal projection. Let us assume that $i=1$.

{\bf The case of four lines.} The class of the free part of the
restriction of the system $\Sigma^2$ onto the surface $S^2$ in
this case is
$$
\xi=nH_S-m_1E_{S,1}-m_2E_{S,2}-\sum_{j=1}^4\mu_jC_j^2.
$$
Since
$H_S=R_1+R_2+\sum_{j=1}^4C_j^2$, $E_{S,1}=R_1$, $E_{S,2}=R_2$,
this class can be re-written in the form
$$
\xi=(n-m_1)R_1+(n-m_2)R_2+\sum_{j=1}^4(n-\mu_j)C_j^2.
$$
Since we are interested in the multiplicities $m_1$ and $\mu_j$,
$j\in \{3,4\}$ only, let us consider in
$$
\Xi={\mathbb Q} R_1\oplus {\mathbb Q} R_2\bigoplus
\mathop{\oplus}\limits_{j=1}^4{\mathbb Q}C_j^2
$$
the orthogonal complement to $C_1^2$, $C_2^2$, $R_2$:
$$
\begin{array}{ccl}
R_1^+ & = & R_1+\frac12 C_1^2+\frac12 C_2^2, \\
C_3^+ & = & C_3^2+\frac12 C_1^2+\frac12 C_2^2+\frac12 R_2,\\
C_4^+ & = & C_4^2+\frac12 C_1^2+\frac12 C_2^2+\frac12 R_2.
\end{array}
$$
The multiplication table for these three classes is given by the
matrix
$$
\frac{1}{2}\left(
\begin{array}{ccc}
-2  &  2   &  2   \\
2   &  -1  &  3   \\
2   &  3   &  -1
\end{array}\right),
$$
the inverse matrix for which is
$$
\frac{1}{12}\left(
\begin{array}{ccc}
-4  &  4   &  4   \\
4   &  -1  &  5   \\
4   &  5   &  -1
\end{array}\right).
$$
Since the sum of the first and the second (or the second and the
third) rows is non-negative, this means precisely that
$$
(n-m_1)+(n-\mu_3)\geqslant 0
$$
(respectively, $(n-m_1)+(n-\mu_4)\geqslant 0$). Proof of
Proposition 2.5 for the case of four lines is complete.

{\bf The case of a conic and two lines.} Here $L=C_3$, and the
multiplication table for the five curves $R_1^2, C_1^2, C_2^2,
C_3^2, R_2^2$ is given by the matrix (we did not give it in
Subsection 4.2):
$$
\left(
\begin{array}{ccccc}
-2 & 1  & 1  & 0  & 0  \\
1  & -2 & 1  & 1  & 1  \\
1  & 1  & -2 & 1  & 0 \\
0  & 1  & 1  & -2 & 1  \\
0  & 1  & 0  & 1  & -2
\end{array}\right).
$$
The orthogonal complement to $C_1^2$, $C_2^2$ is given by the
formulas
$$
\begin{array}{ccl}
R_1^+ & = & R_1^2+C_1^2+ C_2^2, \\
C_3^+ & = & C_3^2+C_1^2+ C_2^2,\\
R_2^+ & = & R_2^2+\frac23 C_1^2+\frac12 C_2^2.
\end{array}
$$
The multiplication table for these classes is
$$
\left(
\begin{array}{ccc}
0  &  2   &  1   \\
2   &  0  &  2   \\
1   &  2   &  -\frac43
\end{array}\right),
$$
and the inverse matrix is
$$
\frac{1}{40}\left(
\begin{array}{ccc}
-12  &  14   &  12   \\
14   &  -3   &  6   \\
12   &  6    &  -12
\end{array}\right).
$$
Since the sum of the first and second rows consists of
non-negative numbers, Proposition 2.5 is shown in this case.

{\bf The case of a cubic curve and a line.} Here the
multiplication table for $R_1^2, C_1^2, C_2^2, R_2^2$ (where $C_1$
is the plane cubic and $C_2$ is the line that we need) is
$$
\left(
\begin{array}{cccc}
-2 & 2  & 0  & 0  \\
2  & -2 & 2  & 1  \\
0  & 2  & -2 & 1  \\
0  & 1  & 1  & -2
\end{array}\right).
$$
The orthogonal complement to $C_1^2$, $R_2^2$ is given by the
formulas
$$
\begin{array}{ccl}
R_1^+ & = & R_1^2+\frac43 C_1^2+ \frac23 R_2^2, \\
C_2^+ & = & C_2^2+\frac53 C_1^2+ \frac43 R_2^2,
\end{array}
$$
the multiplication table for $R_1^+$, $C_2^+$ is
$$
\frac13\left(
\begin{array}{cc}
2    &  10   \\
10   &  8
\end{array}\right),
$$
and the inverse matrix is
$$
\frac{1}{14}\left(
\begin{array}{cc}
-4  &  5 \\
5   &  -1
\end{array}\right).
$$
Q.E.D. for Proposition 2.5.\vspace{0.3cm}


{\bf 4.4. Lines through the point $o_2$.} Let us show Proposition
0.1. In the coordinate notations of Subsection 0.2, using a
suitable linear change of coordinates, we assume that
$o_2=(0,0,1,0)$. The form $q_2(z_1,z_2,z_3)=z_1^2+z_2^2+z_3^2$ is
fixed, $q_2(o_2)\neq 0$, since the line $[o_1,o_2]$ can not be
contained in the tangent cone $T_{o_1}V$ (otherwise,
$[o_1,o_2]\subset V$, contrary to our assumption). Let us find
out, what conditions for the polynomials $q_3(z_*)$ and $q_4(z_*)$
are imposed by the requirement $o_2\in\mathop{\rm Sing} V$. We
will need some routine coordinate computations, which we will
reduce to the minimum, only presenting the results of each step.
Write down
$$
q_3(z_*)=\sum_{i=0}^3 z_3^{3-i}q_{3,i},\quad
q_4(z_*)=\sum_{i=0}^4 z_3^{4-i}q_{4,i},
$$
where $q_{3,i}$ and $q_{4,i}$ are homogeneous polynomials of
degree $i$ in $z_1$, $z_2$, $z_4$. Since the blow-up of the point
$o_1$ resolves that singularity, $q_{3,3}(z_1,z_2,z_4)$ contains
$z_4^3$ with a non-zero coefficient. Consider the new coordinate
system $(z_1,z_2,u,z_4)$, where $z_3=1+u$.

{\bf Lemma 4.1.} {\it The condition $o_2\in V$ is equivalent to
the equality} $q_{3,0}+q_{4,0}=-1$.

{\bf Proof:} obvious computations.

Assuming the equality of the lemma, let us consider the condition
$o_2\in\mathop{\rm Sing} V$.

{\bf Lemma 4.2.} {\it The point $o_2$ is a singularity of the
hypersurface $V$ if and only if the the following equalities
hold:}
$$
3q_{3,0}+4q_{4,0}=-2,\quad q_{3,1}+q_{4,1}\equiv 0.
$$

{\bf Proof:} obvious computations.

From these lemmas we get the equalities $q_{3,0}=-2$, $q_{4,0}=1$.
We also assume that the second equality of Lemma 4.2 holds.

{\bf Lemma 4.3.} {\it In the coordinates $(z_1,z_2,u,z_4)$ the
quadratic part of the polynomial $f$ at the point $o_2$ is
\begin{equation}\label{08.11.2025.1}
w_2=u^2+uq_{4,1}+(z_1^2+z_2^2+q_{3,2}+q_{4,2}),
\end{equation}
the cubic part of the polynomial $f$ at the point $o_2$ is
$$
w_3=2u^3+2u^2q_{4,1}+u(q_{3,2}+2q_{4,2})+(q_{3,3}+q_{4,3}),
$$
and the component of degree 4 is}
$$
q_4(z_1,z_2,u,z_4)=\sum_{i=0}^4 u^{4-i}q_{4,i}(z_1,z_2,z_4).
$$

{\bf Proof:} obvious computations.

It is clear from the explicit formula (\ref{08.11.2025.1}), that
for a general quartic of the second type the rank of the point
$o_2$ is equal to 4. Identifying
$$
{\mathbb E}_2={\mathbb P}_{(z_1:z_2:u:z_4)}^3,
$$
we see that the quadric $Q_2$ is the set of zeros of the
polynomial (\ref{08.11.2025.1}), where the point $p=(0:0:1:0)$
does not lie on $Q_2$. Furthermore, for the fixed forms $q_{4,1}$,
$q_{3,2}$, $q_{4,2}$ the polynomial $w_3$ contains an arbitrary
cubic polynomial $(q_{3,3}+q_{4,3})$ in $(z_1,z_2,z_4)$, that is,
the equation of an arbitrary cubic cone with the vertex at the
point $p$, which implies that for a general quartic of the second
type the equation $w_3|_{Q_2}=0$ defines on the quadric a
non-singular curve of bi-degree (3,3). Finally the polynomial
$q_4(z_1,z_2,u,z_4)$ contains an arbitrary component
$q_{4,4}(z_1,z_2,z_4)$ (that is, the equation of an arbitrary cone
of degree 4 with the vertex at the point $p$), which is not
contained in the equation, defining $Q_2$ and the curve
$\{w_3|_{Q_2}=0\}$. It follows that the system of equations
$$
w_3|_{Q_2}=0,\quad q_4(z_1,z_2,u,z_4)|_{Q_2}=0
$$
has precisely 24 distinct solutions. Q.E.D. for Proposition 0.1.


\section{Infinitely near lines}

In this section we consider only quartics of the first type. We
show Proposition 2.4.\vspace{0.3cm}

{\bf 5.1. The technique of counting multiplicities.} Assume
(contrary to the claim of Proposition 2.4) that the centre of the
maximal singularity $E$ on $V_1$ is a generator $R$ of the cone
$Q_1$. By Proposition 3.1 the singularity $E$ is not the
exceptional divisor of the blow-up of $R$. Let us re-denote the
strict transform $R^2$ by the symbol $B_1$, it is a non-singular
rational curve on the non-singular three-dimensional variety
$V_2$. Let
$$
\sigma_i\colon V^{(i)}\to V^{(i-1)}
$$
be the sequence of blow-ups of curves $B_{i-1}\subset V^{(i-1)}$
with the exceptional divisors $E^{(i)}\subset V^{(i)}$,
$E^{(i)}=\sigma_i^{-1}(B_{i-1})$, where $i=2,\dots, k+1$,
$k\geqslant 2$, such that:

--- $V^{(1)}=V_2$ and $B_1=R^2$,

--- the centre $B_{i-1}$ of each of these blow-ups is the
centre of the maximal singularity $E$ on $V^{(i-1)}$,

--- $E^{(k+1)}$ realizes $E$ on $V^{(k+1)}$ (in other words,
$E^{(k+1)}$ is the centre of $E$ on $V^{(k+1)}$).

Thus this sequence of blow-ups is the resolution of the maximal
singularity $E$ in the sense of \cite[Chapter 2]{Pukh13a}. The
inequality $k\geqslant 2$ holds by Proposition 3.1. Let
$\Sigma^{(i)}$ be the strict transform of the system $\Sigma$ on
$V^{(i)}$. Set
$$
\nu_i=\mathop{\rm mult}\nolimits_{B_{i-1}}\Sigma^{(i-1)},
$$
$i=2,\dots, k+1$, then the Noether-Fano inequality for $E$ (that
is, the condition that $E$ is a maximal singularity of the system
$\Sigma$) takes the form
$$
p_1m_1+\sum_{i=2}^{k+1}p_i\nu_i>n\left( p_1+\sum_{i=2}^{k+1}p_i\right),
$$
where $p_i$ is the number of paths in the graph $\Gamma$ of the
singularity $E=E^{(k+1)}$ from the vertex $(k+1)$ to the vertex
$i$, see \cite[Chapter 2]{Pukh13a} and \cite{Pukh00c}. We will not
dwell on the details, because in the case under consideration this
formalism simplifies essentially: one may assume that
$\nu_{k+1}>n$, whence by the inequality $m_1\leqslant n$ it
follows that the graph $\Gamma$ is a chain, all $p_i$ are equal to
1, for $i\leqslant k$ the curve $B_i\subset E^{(i)}$ is a section
of the projection $E^{(i)}\to B_{i-1}$, so that all curves $B_i$
and all varieties $V^{(i)}$ are non-singular. The Noether-Fano
inequality takes the simplified form
$$
m_1+\nu_2+\dots+\nu_{k+1}>(k+1)n.
$$
Let $Z=(D_1\circ D_2)$ be the self-intersection of the system
$\Sigma$, where $D_1,D_2\in \Sigma$ are general divisors. This is
and effective 1-cycle on $V$, for which, according to the
technique of counting multiplicities, see \cite[Chapter 2,
Proposition 2.4, Corollary 2.5]{Pukh13a}, and the computation at
the end of \cite[Chapter 2, Subsection 2.2]{Pukh13a}, and also
\cite[\S 3]{Pukh2025a}, we have the estimate
$$
\mathop{\rm mult}\nolimits_{o_1} Z>\frac{2(k+1)^2}{2k+1}n^2.
$$
It is easy to check that the right hand side of this inequality is
strictly higher than $4n^2$ for $k\geqslant 3$. Since $\mathop{\rm
deg} Z=4n^2$, we conclude that $k=2$. Thus the maximal singularity
is realized by the blow-up of the curve $B_2$ (a section of the
ruled surface $E^{(2)}\to B_1=R^2$), and the Noether-Fano
inequality takes the form of the estimate
\begin{equation}\label{13.11.2025.1}
m_1+\nu_2+\nu_3>3n.
\end{equation}
However, for this (only this) case $k=2$ the standard technique of
counting multiplicities does not make it possible to exclude the
maximal singularity, and we need to study this case in more
detail.\vspace{0.3cm}


{\bf 5.2. Curves on the quadric $Q_2$.} Let $D\in\Sigma$ be a
general divisor, $D^2$ its strict transform on $V_2$. Since
$D^2\sim nH-m_1Q_1-m_2Q_2$, the restriction $D^2|_{Q_2}$ is an
effective 1-cycle on $Q_2\cong {\mathbb P}^1\times {\mathbb P}^1$
of bi-degree $(m_2,m_2)$ (in other words, this 1-cycle is cut out
on $Q_2$ by a surface of degree $m_2$ in ${\mathbb E}_2\cong
{\mathbb P}^3$).

{\bf Proposition 5.1.} {\it The following estimate is true:}
\begin{equation}\label{13.11.2025.2}
\nu_2+\nu_3\leqslant 3m_2.
\end{equation}

{\bf Proof.} Set
$$
p_1=B_1\cap Q_2
$$
(the intersection is transversal, since $B_1=R^2$ intersects
${\mathbb E}_2$ transversally at this point) and let
$N=\sigma_2^{-1}(p_1)\subset E^{(2)}$ be the fibre of the ruled
surface $E^{(2)}$ over the point $p_1\in B_1$. Obviously,
$Q_2^{(2)}$ (the strict transform of $Q_2$ on $V^{(2)}$) is the
result of blowing up the point $p_1\in Q_2$, and $N=Q_2^{(2)}\cap
E^{(2)}$ is the exceptional line of that blow-up. Furthermore, the
restriction of the divisor
$$
D^{(2)}\sim nH-m_1Q_1-m_2Q_2-\nu_2E^{(2)}
$$
onto the surface $Q_2^{(2)}$ is precisely
$$
\left(\sigma_2|_{Q_2^{(2)}}\right)^*D^2|_{Q_2}-\nu_2 N.
$$
This effective 1-cycle on $Q_2^{(2)}$ has multiplicity $\geqslant
\nu_3$ at the point $p_2=B_2\cap Q_2^{(2)}=B_2\cap N$ (the
intersections are transversal, since $B_2$ is a section of the
ruled surface $E^{(2)}/B_1$).

Therefore, we obtain the following geometric situation. On the
non-singular quadric $Q_2\subset {\mathbb E}_2$ there is an
effective 1-cycle $Y$ of bi-degree $(m_2,m_2)$, satisfying the
inequality $\mathop{\rm mult}\nolimits_{p_1} Y\geqslant \nu_2$. On
the blow-up of the point $p_1\in Q_2$ the effective 1-cycle
$Y-\nu_2 N$ (recall that $N$ is the exceptional line of that
blow-up, and the pull back of $Y$, as a divisor, from $Q_2$ onto
that blow-up is denoted by the same symbol) has multiplicity
$\geqslant \nu_3$ at the point $p_2\in N$. We will consider the
point $p_2$ as a point, infinitely near to $p_1$, in the sense of
the blow-up of the point $p_1$ on ${\mathbb E}_2\cong {\mathbb
P}^3$, and let $[p_1,p_2]\subset {\mathbb E}_2$ be the line on
${\mathbb E}_2$, passing through $p_1$ in the direction of $p_2$.
There are two options:

(1) $[p_1,p_2]\not\subset Q_2$, and then $[p_1,p_2]\cap
Q_2=\{p_1\}$ (set-theoretically),

(2) $[p_1,p_2]\subset Q_2$ is one of the two lines on $Q_2$,
passing through the point $p_1$.

Assume that the case (1) takes place. Let $\Lambda\subset Q_2$ be
the section of $Q_2$ by a general plane in ${\mathbb E}_2$,
containing the line $[p_1,p_2]$. Then the intersection of
$\Lambda$ with the support of the cycle $Y$ is zero-dimensional,
so that we get:
$$
2m_2=(\Lambda\cdot Y)\geqslant (\Lambda\cdot Y)_{p_1}\geqslant
\nu_2+\nu_3,
$$
which is considerably stronger than the inequality of Proposition
5.1.

Assume that the case (2) takes place. Write down
$Y=a[p_1,p_2]+Y^+$, where $a\in {\mathbb Z}_+$ and the effective
1-cycle $Y^+$ of bi-degree $(m_2-a,m_2)$ does not contain the line
$[p_1,p_2]$ as a component. We get
$$
m_2=([p_1,p_2]\cdot Y^+)\geqslant \mathop{\rm mult}\nolimits_{p_1} Y^++
\mathop{\rm mult}\nolimits_{p_2}(Y^+)^{(2)},
$$
where the symbol $(Y^+)^{(2)}$ means the strict transform of the
1-cycle $Y^+$ on the blow-up of the point $p_1$ on $Q_2$ (or on
$V^{(2)}$, which is the same). Now we get:
$$
\nu_2+\nu_3\leqslant \mathop{\rm mult}\nolimits_{p_1} Y+
\mathop{\rm mult}\nolimits_{p_2} Y^{(2)}=
\mathop{\rm mult}\nolimits_{p_1} Y^++
\mathop{\rm mult}\nolimits_{p_2}(Y^+)^{(2)}+2a\leqslant 3m_2,
$$
where $Y^{(2)}$ means the corresponding strict transform. This
completes the proof of Proposition 5.1.\vspace{0.3cm}


{\bf 5.3. The two techniques of the proof.} For a quartic of the
first type of general position there are only two options, see
Subsections 3.1 and 3.5:

(1) the case of general position, when $(V\circ P)=V\cap P$ is an
irreducible curve $C$ of degree 4 with the triple point $o_1$, so
that $o_2\in C^1$ is a non-singular point (this case takes place
for almost all generators $R$, except for finitely many of them),

(2) $(V\circ P)=C+L$, where $C$ is an irreducible cubic curve with
the double point $o_1$, $L\ni o_1$  is a line, and moreover,
$o_2\not\in L^1$ and $o_2\in C^1$ is a non-singular point.

We will exclude these two cases in different ways: the method of
excluding the case (1) does not work in the case (2), and the
method of excluding the case (2) does not apply to finitely many
generators $R$ and planes $P$, for which the curve $C$ has at the
point $o_1$ two branches, a non-singular one and a cuspidal one,
and moreover, the point $o_2$ lies exactly on the cuspidal branch.
The case (1) will be excluded by means of the ``quadratic''
technique, which is more traditional and simpler. Our approach to
the proof of Proposition 2.4 in the case (1) is to estimate the
input of the curve $C$ to the self-intersection $Z$ of the system
$\Sigma$; we will show that this input is sufficiently big, which
will make it possible to obtain a contradiction by means of the
standard technique of counting multiplicities:
\begin{equation}\label{13.11.2025.3}
\mathop{\rm mult}\nolimits_{o_1} Z\geqslant 2m_1^2+\nu_2^2+\nu_3^2,
\end{equation}
see \cite[\S 3]{Pukh2025a}, although it does not give the
$4n^2$-inequality.

In the case (2) (that is, for finitely many generators $R$ of the
cone $Q_1$) this approach does not give a contradiction, and we
will argue in a different way: using the linear technique of \S 3,
we will estimate the multiplicity of the curve $C$ with respect to
the system $\Sigma$ and make the estimate of Proposition 5.1
stronger, which will lead to a contradiction without the technique
of counting multiplicities.

We start with the case of general position.\vspace{0.3cm}


{\bf 5.4. The case of an irreducible plane quartic.} We consider
the situation of Subsection 3.5. In the notations of Subsection
3.5 we have on the non-singular surface $S^2$ (where $S\in |H-P|$
is a general element of the pencil):
$$
H_S=C^2+2R_1^2+R_2^2+2R_3
$$
is the pull back of the hyperplane section,
$$
E_{1,S}=R_1^2+R_2^2+R_3,\quad E_{2,S}=R_3.
$$
The class of the mobile part of the restriction of the system
$\Sigma^2$ onto $S^2$ is
$$
\xi=nH_S-m_1E_{1,S}-m_2E_{2,S}-bR_1^2-b_C C^2,
$$
where $b_C=\mathop{\rm mult}\nolimits_C\Sigma$, and the other
coefficients have the same meaning as in \S 3; of course,
$b=\nu_2$. Re-writing
$$
\xi=(n-b_C)C^2+(2n-m_1-b)R_1^2+(n-m_1)R_2^2+(2n-m_1-m_2)R_3
$$
and multiplying by $C^2$, we get (see Subsection 3.5):
\begin{equation}\label{13.11.2025.4}
b_C\geqslant \frac12(3m_1+m_2+2\nu_2-4n).
\end{equation}
(It follows from the inequalities (\ref{13.11.2025.1}) and
(\ref{13.11.2025.2}) that the right hand side is strictly
positive.) It is easy to check that for a fixed value
$\nu_2+\nu_3=2\theta$ the estimates (\ref{13.11.2025.3}) and
(\ref{13.11.2025.4}) are the worst for $\nu_2=\nu_3=\theta$
(recall that $\nu_2\geqslant \nu_3$: on a non-singular variety
multiplicities do not increase), so that we may assume that this
equality takes place.

Consider first the case $\theta\in (n,\frac98 n]$. In that case
the inequality
$$
3n-2\theta\geqslant \frac23\theta
$$
holds, so that, making the estimates (\ref{13.11.2025.3}) and
(\ref{13.11.2025.4}) worse, set $m_1=3n-2\theta$ and $m_2=\frac23
\theta$, whence we get:
$$
b_C\geqslant \frac16 (15n-10\theta)
$$
and
$$
\mathop{\rm mult}\nolimits_{o_1} Z\geqslant 18n^2-24n\theta+
10\theta^2.
$$
Write down $Z=\beta C+Z^{\sharp}$, where the effective 1-cycle
$Z^{\sharp}$ does not contain the curve $C$ as a component and
$\beta\geqslant b_C^2$. Obviously,
$$
\mathop{\rm deg} Z=4\beta+\mathop{\rm deg} Z^{\sharp}=4n^2
$$
and
$$
\mathop{\rm mult}\nolimits_{o_1} Z\leqslant 3\beta+\mathop{\rm
deg} Z^{\sharp}= 4n^2-\beta.
$$
Combining these estimates, we obtain the inequality
$$
18n^2-24n\theta+10\theta^2\leqslant 4n^2-
\frac{1}{36}(15n-10\theta)^2,
$$
which after easy transformations takes the form
$460\theta^2-1164n\theta+729 n^2\leqslant 0$, where
$\theta\in[n,\frac98 n]$. However, the quadratic function
$460t^2-1164t+729$ is decreasing on the interval $[1,\frac98]$ and
its value at the right hand end of this interval is
$\frac{108}{64}>0$. This contradiction excludes the case
$\theta\leqslant \frac98 n$.

Assume now that $\theta>\frac98 n$. In that case
$3n-2\theta<\frac23\theta$, so that, making the estimates
(\ref{13.11.2025.3}) and (\ref{13.11.2025.4}) worse, we set
$m_1=m_2=\frac23\theta$. From there it follows that
$$
b_C\geqslant \frac13(7\theta-6n)
$$
and $\mathop{\rm
mult}\nolimits_{o_1}Z\geqslant\frac{26}{9}\theta^2$. Again we
write down $Z=\beta C+Z^{\sharp}$ with $\beta\geqslant b_C^2$ and
obtain (in the same way as above) the inequality
$139\theta^2-300n\theta+162n^2\leqslant 0$, where
$\theta\geqslant\frac98 n$. However the quadratic function
$139t^2-300t+162$ is strictly increasing for $t\geqslant\frac98$,
and its value for $t=\frac98$ is $\frac{27}{64}>0$ (up to a
factor, it is the same value as above, since for $\theta=\frac98n$
we have $\frac23\theta=3n-2\theta$). Contradiction.

Proof of Proposition 2.4 in the case of an irreducible plane
quartic is complete.\vspace{0.3cm}


{\bf 5.5. The case of a plane cubic and a line.} In the case (2)
we use the constructions and notations of Subsection 3.2, in
certain situations modifying the notations, which will be
specially noted. The point $p_1\in Q_2$ and the infinitely near
(to $p_1$) point $p_2$ have the same meaning as in Subsection 5.2.
The curve $C^1$ is non-singular. Since the case (2) corresponds to
finitely many generators $R$ (and finitely many planes $P$), we
can say more about the curve $C$: it is a cubic curve with a nodal
singularity $o_1$, that is, it has two branches at this point,
which meet transversally at the point $o_1$. This is true for all
generators $R$ in the case (2) for a (Zariski) general quartic of
the first type. Therefore, the strict transform $C^2$ intersects
transversally the exceptional divisor ${\mathbb E}_2$ and the
quadric $Q_2$ at the point $p_3$, where $p_3\neq p_1$. Since
$(P^2\circ V_2)$ does not contain the exceptional line $P^2\cap
{\mathbb E}_2$ (because the equality $\mathop{\rm
mult}\nolimits_{o_2} (P^1\circ V^1)=2$ holds), we conclude that
the line
$$
[p_1,p_3]=P^2\cap {\mathbb E}_2,
$$
containing the points $p_1,p_3$ in ${\mathbb E}_2\cong {\mathbb
P}^3$, is not contained in the non-singular quadric $Q_2$. It
follows that $p_3\not\in T_{p_1}Q_2$. Set $b_C=\mathop{\rm
mult}\nolimits_C\Sigma$. Now we can state a stronger version of
Proposition 5.1.

{\bf Proposition 5.2.} {\it The following estimate holds:}
\begin{equation}\label{15.11.2025.1}
\nu_2+\nu_3+b_C\leqslant 3m_2.
\end{equation}

{\bf Proof.} We use the notations of Proposition 5.1. Obviously,
$\mathop{\rm mult}\nolimits_{p_3} Y=b^*\geqslant b_C$. If
$[p_1,p_2]\not\subset Q_2$, then we consider the section $\Lambda$
of the quadric $Q_2$ by a plane, containing the line $[p_1,p_2]$
and the point $p_3$. This plane is unique precisely because
$p_3\not\in [p_1,p_2]$. The conic $\Lambda$ is certainly
non-singular. Write down
$$
Y=l\Lambda+Y^{\sharp},
$$
where the 1-cycle $Y^{\sharp}$ is effective and does not contain
$\Lambda$ as a component. We have:
$$
2(m_2-l)=(\Lambda\cdot Y^{\sharp})\geqslant
(\Lambda\cdot Y^{\sharp})_{p_1}+
(\Lambda\cdot Y^{\sharp})_{p_3}\geqslant (\nu_2-l)+(\nu_3-l)+(b_C-l),
$$
so that
$$
\nu_2+\nu_3+b_C\leqslant 2m_2+l\leqslant 3m_2.
$$
Assume that $[p_1,p_2]\subset Q_2$. Set $L_1=[p_1,p_2]$ and denote
by the symbol $L_2$ the only line in the transversal to $L_1$
family of lines on $Q_2$, passing through the point $p_3$ (that
is, if $L_1$ is of bi-degree $(1,0)$, then $L_2$ is of bi-degree
$(0,1)$). Write down
$$
Y=l_1L_1+l_2L_2+Y^{\sharp},
$$
where the effective 1-cycle $Y^{\sharp}$ does not contain $L_1$,
$L_2$ as components and $l_1,l_2\in {\mathbb Z}_+$. Obviously,
$l_1\leqslant m_2$ and $l_2\leqslant m_2$ and $Y^{\sharp}$ is of
bi-degree $(m_2-l_1,m_2-l_2)$. Now we get:
$$
2m_2-l_1-l_2=((L_1+L_2)\cdot Y^{\sharp})\geqslant
\mathop{\rm mult}\nolimits_{p_1} Y^{\sharp}+
\mathop{\rm mult}\nolimits_{p_2} (Y^{\sharp})^{(2)}+
\mathop{\rm mult}\nolimits_{p_3} Y^{\sharp}\geqslant
$$
$$
\geqslant \nu_2+\nu_3+b_C-2l_1-l_2,
$$
whence we obtain
$$
\nu_2+\nu_3+b_C\leqslant 2m_2+l_1\leqslant 3m_2.
$$
Q.E.D. for the proposition.

Set $b_L=\mathop{\rm mult}\nolimits_{L}\Sigma$. Replacing in the
notations of Subsection 3.2 the multiplicity $b$ by $\nu_2$, the
curve $C_1$ by $C$ and the line $C_2$ by $L$, we get the effective
class
$$
\xi=nH_S-m_1E_{1,S}-m_2E_{2,S}-\nu_2R^2_1-b_CC^2-b_LL^2=
$$
$$
=(n-b_C)C^2+(n-b_L)L^2+(2n-m_1-\nu_2)R^2_1+(n-m_1)R^2_2+(2n-m_1-m_2)R_3
$$
on the non-singular surface $S^2$. The inequalities $(\xi\cdot
C^2)\geqslant 0$ and $(\xi\cdot L^2)\geqslant 0$ turn into the
estimates (see the multiplication table (\ref{03.11.2025.1}))
$$
b_C\geqslant\frac12(b_L+2m_1+m_2+\nu_2-3n),\quad
b_L\geqslant\frac12(b_C+m_1+\nu_2-n).
$$
Substituting the right hand side of the second inequality for
$b_L$ in the right hand side of the first one, after easy
computations we get
\begin{equation}\label{17.11.2025.1}
b_C\geqslant\frac13(5m_1+2m_2+3\nu_2-7n).
\end{equation}
Besides, Proposition 5.2 gives the estimates
$$
m_1\geqslant m_2\geqslant\frac13(\nu_2+\nu_3+b_C)
$$
Substituting $\frac13(\nu_2+\nu_3+b_C)$ for $m_1,m_2$ in
(\ref{17.11.2025.1}), we get
$$
2b_C\geqslant 16\nu_2+7\nu_3-21n>2n
$$
(since $\nu_2\geqslant\nu_3>n$), so that $b_C>n$, which
contradicts the plane section lemma (Proposition 2.2).

Proof of Proposition 2.4 is complete.

Q.E.D. for Theorems 0.2 and 0.3.


\begin{flushleft}
Department of Mathematical Sciences,\\
The University of Liverpool
\end{flushleft}

\noindent{\it pukh@liverpool.ac.uk}

\end{document}